%V.A. Gritsenko and V.V. Nikulin
%Automorphic forms and Lorentizian Kac--Moody algebras. Part I. 
%till March 1997: gritsenk@kurims.kyoto-u.ac.jp  nikulin@kurims.kyoto-u.ac.jp

\input amstex
\documentstyle{amsppt}
\magnification=1200
\catcode`\@=11
\redefine\logo@{}
\catcode`\@=13

\define \bn{\Bbb N}
\define \bz{\Bbb Z}
\define \bq{\Bbb Q}
\define \br{\Bbb R}
\define \bc{\Bbb C}

\define \M{{\Cal M}}
\define\Ha{{\Cal H}}
\define\La{{\Cal L}}
\define\geg{{\goth g}}
\define\0o{{\overline 0}}
\define\1o{{\overline 1}}
\define\rk{\text{rk}~}

%additional to the file

\define\io{{\overline i}}
\define\mult{\text{mult}}
\define\re{{\text{re}}}
\define\im{{\text{im}}}

%additional, sept 14, 1996

\define\Sym{{\text{Sym}~}}

\input epsf.tex

\TagsOnRight

\document
\document

\topmatter
\title
Automorphic Forms and Lorentzian Kac--Moody algebras. Part I.  
\endtitle

\author
Valeri A. Gritsenko \footnote{Supported by
RIMS of Kyoto University 
\hfill\hfill}
and
Viacheslav V. Nikulin \footnote{Supported by
Grant of Russian Fund of Fundamental Research and RIMS of Kyoto University 
\hfill\hfill}
\endauthor

\address
St. Petersburg Department of Steklov Mathematical Institute,
\newline
${}\hskip 8pt $
Fontanka 27, 191011 St. Petersburg,  Russia
\endaddress
\email
gritsenk\@gauss.pdmi.ras.ru 
\endemail

\address
Steklov Mathematical Institute,
ul. Vavilova 42, Moscow 117966, GSP-1, Russia
\endaddress

\email
slava\@nikulin.mian.su
\endemail

\abstract
Using the general method which was applied to prove 
finiteness of the set of hyperbolic generalized Cartan matrices 
of elliptic and parabolic type, we classify all symmetric 
(and twisted to symmetric) hyperbolic generalized Cartan matrices 
of elliptic type and of rank $3$ with a lattice Weyl vector. 

We develop the general theory of reflective lattices $T$ with 
$2$ negative squares and reflective automorphic forms on 
homogeneous domains of type IV defined by $T$. We consider this 
theory as mirror symmetric to the theory of elliptic and 
parabolic hyperbolic reflection groups and corresponding 
hyperbolic root systems. We formulate Arithmetic Mirror Symmetry 
Conjecture relating both these theories and prove some  
statements to support this Conjecture. This subject is connected 
with automorphic correction of Lorentzian Kac--Moody algebras. 
We define Lie reflective automorphic forms which are the 
most beautiful automorphic forms defining automorphic Lorentzian 
Kac--Moody algebras and formulate finiteness Conjecture for these forms. 

Detailed study of automorphic correction and Lie reflective automorphic 
forms for generalized Cartan matrices mentioned above will be given 
in Part II.  
\endabstract

\rightheadtext
{Lorentzian Kac--Moody algebras}
\leftheadtext{V. Gritsenko and  V. Nikulin}
\endtopmatter

\document

\head
0. Introduction
\endhead

This paper continues series of our papers 
\cite{GN1}---\cite{GN4}, \cite{N9}---\cite{N11} devoted to  
automorphic forms related with   
Lorentzian Kac--Moody algebras and Mirror Symmetry. 
 
In this paper we try to develop regularly the theory of Lorentzian 
Kac--Moody algebras and the related theory of automorphic forms. 

Like the classical theory of semi-simple Lie algebras and affine Kac--Moody 
algebras, the theory of Lorentzian Kac--Moody 
algebras and the related theory of 
automorphic forms start with some root systems and reflection groups. 
For Lorentzian Kac--Moody algebras, they are hyperbolic reflection groups 
and corresponding root systems which are given by 
reflection groups of hyperbolic lattices (over $\bz$). It was shown in 
\cite{N9}, \cite{N10} that to have a Lorentzian Kac--Moody algebra with 
good properties one is forced to consider hyperbolic reflection groups  
having some condition of finiteness of volume for a fundamental 
polyhedron of the group. Hyperbolic 
reflection groups of this type (and corresponding hyperbolic 
root systems) are divided in 
{\it elliptic,} when volume is finite, and 
{\it parabolic,} when volume is finite after 
factorization of a fundamental polyhedron by its 
Abelian cusp group of symmetries. Hyperbolic lattices having 
these groups as reflection groups are respectively called {\it reflective}  
(elliptic reflective or parabolic reflective). It was shown 
in \cite{N9} and \cite{N10} (in fact, it was done in \cite{N4}, \cite{N5} 
long time ago)  
that sets of elliptic and parabolic hyperbolic reflection groups and 
corresponding hyperbolic generalized Cartan matrices of elliptic and 
parabolic type  
(which define the corresponding Lorentzian Kac--Moody algebras) 
are essentially finite. More exactly, the set of reflective hyperbolic 
lattices is finite for $\rk \ge 3$. It shows that ``right'' 
Lorentzian Kac--Moody algebras are 
all {\it exceptional}, there are no infinite series (at least, for 
elliptic type).  

The main idea of proving finiteness of elliptic and parabolic 
hyperbolic reflection groups is not only theoretical. This gives 
an effective method for classification.  It has been
applied in \cite{N7} to classify all hyperbolic lattices having
elliptic groups of reflections in vectors with equal squares (so
called elliptic 2-reflective hyperbolic lattices).  
(For rank $\ge 5$ this classification had been 
obtained in \cite{N3}; for rank $=4$ 
it had been done by \'E.B. Vinberg, see \cite{N8}.)  
We mention that
this gives a description of all symmetric hyperbolic generalized Cartan
matrices of elliptic type and of rank $3$. They are divided in a
finite set of series (possibly infinite) corresponding to their 
2-reflective hyperbolic lattices.  
These hyperbolic lattices are important for the
theory of K3 surfaces (e.g., for automorphism groups and for 
discriminants of moduli).  
For automorphic correction of Lorentzian
Kac--Moody algebras, which we discuss below, one needs to find among 
these matrices all matrices with so called lattice Weyl vector
(possibly for a twisted generalized Cartan matrix). Only for these
matrices one may hope to find an automorphic form on IV type
domains which gives an automorphic correction of the corresponding
Kac--Moody algebra.  Here we don't discuss possible non-holomorphic
automorphic forms, forms on different type domains or possible
vector-valued holomorphic forms. From our point of view, the case of
lattice Weyl vector is the most symmetric and beautiful to be
considered.

In Sect. 1, we apply this general method to classify all symmetric
(and twisted to symmetric) hyperbolic generalized Cartan matrices of
elliptic type  and of rank $3$ with a lattice Weyl vector. 
One needs existence of a lattice Weyl vector to 
have a formal denominator function
of the corresponding Kac--Moody algebra as a function on IV type
domain. We prove that there are exactly $16$ symmetric hyperbolic
generalized Cartan matrices of elliptic type with a lattice Weyl
vector and of the rank 3.  Between them $12$ are non-compact (i.e. a
fundamental polyhedron of the corresponding hyperbolic reflection 
group has a vertex at infinity). The list of these $12$ non-compact 
cases was announced in \cite{GN4}. We also find conjecturally all
hyperbolic generalized Cartan matrices of elliptic type and of 
rank $3$ with a lattice Weyl vector which are twisted to symmetric 
generalized Cartan matrices. There are 60 of them (including 16 above). 
All these cases are important because of their relation with
discriminants of K3 surfaces moduli which we will discuss in Part
2 of the paper. We classify corresponding parabolic cases in  
Part 2.

We are sure that similar method permits to find all hyperbolic
generalized Cartan matrices of elliptic and parabolic 
type and of rank 3 with a lattice Weyl vector, 
and there should not be too many of them.  
One can adapt this method for higher-dimensional case
(with combination of other methods). It is expected that number of
possibilities is decreasing when dimension is increasing. Thus, our
calculations suggest that it is reasonable to find all hyperbolic
generalized Cartan matrices of elliptic and parabolic type with a
lattice Weyl vector. There should not be too many cases.

\smallpagebreak 

In Sect. 2 we discuss the generalization of the theory of reflective
hyperbolic lattices and corresponding elliptic and parabolic 
hyperbolic reflection groups and
their elliptic and parabolic hyperbolic 
root systems on complex Hermitian
symmetric domains of type IV. Using Mirror Symmetry ideology (related, 
at least, with K3 surfaces), we are trying to find the corresponding
mirror symmetric theory. 

From our point of view, this mirror symmetric theory is related with
so called {\it reflective automorphic forms} $\Phi$ on IV type domains
$\Omega (T)$ defined by lattices $T$ with $2$ negative squares, i.e. 
of signature $(n, 2)$.

We call an automorphic form $\Phi$ on $\Omega (T)$ with respect
to a subgroup $G\subset O^+(T)$ of finite index {\it reflective} if it
is a holomorphic automorphic form 
and the set of zeros of $\Phi$ is union of quadratic
divisors orthogonal to primitive roots of $T$. Here a primitive element
$\delta \in T$ is called a {\it root} if $(\delta, \delta)>0$ and the
reflection $s_\delta$ in $\delta$ belongs to the automorphism group of
the lattice $T$. A lattice $T$ with $2$ negative squares having a
reflective automorphic form $\Phi$ is called {\it reflective} also.
We consider the set $\Delta (\Phi)$ of all primitive roots of 
$T$ orthogonal to
components of the zero divisor of $\Phi$ as a
generalization of hyperbolic primitive root systems on lattices with
$2$ negative squares. The geometry of the automorphic 
form $\Phi$ should define a
more detailed root system with the primitive root system $\Delta (\Phi
)$ as the base.

We propose {\it Arithmetic Mirror Symmetry Conjecture} which relates
reflective hyperbolic lattices and elliptic and parabolic 
hyperbolic root systems with reflective
lattices $T$ with $2$ negative squares 
and primitive root systems $\Delta (\Phi )$ of their 
reflective automorphic forms $\Phi$. 
(This Conjecture is a generalization of our Mirror Symmetry
Conjectures and results for K3 in
\cite{GN1}---\cite{GN4} and \cite{N10}, \cite{N11}.)  
Roughly speaking, Arithmetic Mirror Symmetry Conjecture says that {\it
a reflective lattice $T$ with $2$ negative squares and a root system
$\Delta (\Phi)$ of a reflective automorphic form $\Phi$ of $T$ should
give hyperbolic reflective lattices and elliptic and parabolic
hyperbolic root systems at cusps (here and in what follows 
we consider only $0$-dimensional cusps). 
All elliptic and parabolic reflective lattices and hyperbolic   
root systems appear at cusps of reflective
lattices $T$ and systems of primitive roots $\Delta (\Phi )$ of their
reflective automorphic forms $\Phi$. 
Number of reflective lattices $T$ of $\rk T\ge 5$ is finite.}

The main reason why Arithmetic Mirror Symmetry Conjecture is valid is
the finiteness results above about reflective hyperbolic lattices and
elliptic and parabolic hyperbolic root systems, and proper formulated
Necessary condition 2.2.2 (see Necessary condition 2.2.2 in Sect. 2)
imposed by the Koecher principle. It was shown in \cite{N11} that
Necessary condition 2.2.2 gives extremely strong restrictions on a
lattice $T$ to have a reflective automorphic form.

We prove two statements which support Arithmetic Mirror Symmetry
Conjecture. First we show that lattices $T$ and root systems
satisfying Arithmetic Mirror Symmetry Conjecture almost satisfy
Necessary condition 2.2.2 (see Theorem 2.2.5). We have mentioned that
this condition is extremely strong. The second argument is related with
Fourier expansions of reflective automorphic forms $\Phi$ at cusps. We
show that if an automorphic form $\Phi$ has a Fourier 
expansion at a cusp with a 
{\it generalized lattice Weyl vector} (see Sect. 2.4), then Arithmetic
Mirror Symmetry Conjecture is valid for $\Phi$ at the cusp. 
This is exactly the same argument which was applied in
\cite{N10}, \cite{N9} 
to show that to have a good Lorentzian Kac--Moody algebra 
(e.g., having an automorphic correction) one is forced to consider
hyperbolic generalized Cartan matrices of elliptic or parabolic 
type with a lattice Weyl vector.  
Infinite products, recently
introduced and studied  by R. Borcherds \cite{B5} and which are 
connected with automorphic forms having rational quadratic divisors, 
give examples of Fourier expansions with generalized lattice Weyl 
vectors. (See also the important paper by J. Harvey and G. Moore 
\cite{HM1} 
closely related with this subject.) Thus, if
$\Phi$ has this infinite product expansion, it satisfies Arithmetic
Mirror Symmetry Conjecture.  Since sets of zeros of reflective
automorphic forms $\Phi$ are unions of quadratic divisors, it is
possible that reflective automorphic forms 
$\Phi$ always have infinite product expansions of this type at cusps. 

We mention that Arithmetic Mirror Symmetry Conjecture together with
finiteness results about elliptic and parabolic hyperbolic reflection
groups and reflective hyperbolic 
lattices are extremely important for finding all
reflective lattices $T$ with $2$ negative squares and their reflective
automorphic forms $\Phi$.

In Sect. 2.5 we consider automorphic forms $\Phi$ having Fourier
expansion at a cusp with integral Fourier 
coefficients and connected with a
hyperbolic generalized Cartan matrix $A$ with a lattice Weyl vector
(by \cite{N10} the matrix $A$ should have elliptic or parabolic type;
otherwise, the automorphic form $\Phi$ does not exist).  We call this
Fourier expansion {\it Fourier expansion of Lie type}. Fourier
expansions of this type were considered in
\cite{N10} and \cite{GN1}---\cite{GN4} where, in particular,  
it was shown that some very classical automorphic forms have Fourier
expansions of this type.  In connection with generalized Kac--Moody
algebras theory, first examples of automorphic forms with Fourier
expansions of Lie type were found by R. Borcherds
(\cite{B4}---\cite{B6}).  A Fourier expansion, at a cusp, of Lie type  
of an automorphic form 
defines generalized Lorentzian Kac--Moody superalgebras which have the
denominator function defined by this expansion.  One can call these
algebras as {\it automorphic Lorentzian Kac--Moody superalgebras}.  Using 
terminology introduced in \cite{N10}, \cite{GN1}---\cite{GN4},
this automorphic form, considered at the cusp, and its 
generalized Lorentzian Kac--Moody superalgebras give {\it
automorphic corrections of the Kac--Moody algebra with the generalized
Cartan matrix $A$}.
 
In Sect. 2.5 we define the most interesting (from our point of view)
class of {\it Lie reflective automorphic forms}. They are Lie
reflective automorphic forms $\Phi$ having Fourier expansions of Lie
type at all cusps.  All automorphic forms with Fourier expansions of
Lie type found in our papers \cite{GN1}---\cite{GN4} are automorphic
forms of this type (Lie reflective).  Main Conjecture is that {\it
there exists only finite number of Lie reflective automorphic forms}.
This Conjecture should not be very difficult because it is very
closely related with finiteness results for hyperbolic generalized
Cartan matrices of elliptic and parabolic type and with a lattice Weyl
vector. Some complete list of these matrices for rank $3$ we have in
Sect. 1.

\smallpagebreak 

We suppose that all finiteness results and conjectures above 
are related with conjectured finiteness of families of 
Calabi--Yau 3-folds and $n$-folds for $n \ge 3$. One can find 
some arguments which support this thesis in \cite{N12}, 
\cite{HM1}, \cite{HM2}, \cite{Ka1}, \cite{Ka2}, \cite{DVV}. 
We hope to consider this point in further publications.  

\smallpagebreak

In Part II of our paper (it will appear in a short time) we find
automorphic corrections of Lorentzian Kac--Moody algebras with the most 
part of 60
hyperbolic generalized Cartan matrices of elliptic type 
classified in Sect. 1 and for some hyperbolic generalized Cartan 
matrices of parabolic type.  For all 
their correcting automorphic forms $\Phi$ we  
prove that they are
Lie reflective finding their Fourier expansions and 
infinite product expansions. 
Using general theory described in Sect. 2, 
we can see that these automorphic forms are not just accidental. 
They are very
universal and natural from the point of view of the general theory of
reflective automorphic forms and elliptic and parabolic 
hyperbolic root systems described above.  From this point of view, 
their Fourier expansions at cusps 
give very natural automorphic Lorentzian Kac--Moody superalgebras. 

This paper was written during our stay at RIMS of Kyoto University. 
We are grateful to the Mathematical Institute for hospitality.

\head
1. Classification of symmetric hyperbolic generalized Cartan matrices
of rank 3 of elliptic arithmetic type and with a lattice Weyl vector
\endhead 

\subhead
1.1.  Reminding of hyperbolic generalized Cartan matrices
\endsubhead 

In this Section, 
we recall elementary facts about connection between hyperbolic
generalized Cartan matrices and reflection groups in hyperbolic
spaces. See \cite{K1} about general definitions of generalized Cartan
matrices and their properties, \cite{V3} about elementary properties
of hyperbolic reflection groups, which we need here, and \cite{N9},
\cite{N10}. 

For a countable set of indices $I$, a finite rank matrix $A=(a_{ij})$,
$i,j \in I$, is called a {\it generalized Cartan matrix} if

\noindent 
(C1)\ \ \ $a_{ii}=2$ for any $i\in I$;

\noindent
(C2) $a_{ij}$ are non-positive integers for $i\not=j$; 
\hskip6cm (1.1.1) 

\noindent
(C3) $a_{ij}=0$ implies $a_{ji}=0$. 

We always suppose that $A$ is {\it indecomposable} which means that
there does not exist a decomposition $I=I_1 \cup I_2$ such that both
$I_1$ and $I_2$ are non-empty and $a_{ij}=0$ for any $i \in I_1$ and
any $j \in I_2$.

A generalized Cartan matrix $A$ is called {\it symmetrizable} if there
exists an invertible diagonal matrix
$D=\text{diag}(...,\epsilon_i,...)$ and a symmetric matrix
$B=(b_{ij})$, such that $$ A=DB;\ \ \text{or}\ \
(a_{ij})=(\epsilon_ib_{ij})
\tag{1.1.2}
$$ One always can suppose that $$
\epsilon_i \in \bq,\ \ \epsilon_i>0,\ \  
b_{ij} \in \bz; \ \ b_{ii} \in 2\bz,
\tag{1.1.3} 
$$ equivalently, $$ b_{ij} \in \bz,\ \ b_{ii} \in 2\bz,\ \ b_{ij}\le
0,\ \ b_{ii}>0,
\tag{1.1.4} 
$$ for any $i,j \in I$. 

Since the matrix $A$ is indecomposable, the matrices $D$ and $B$ are
defined uniquely up to a multiplicative constant.  
We call the matrix $B$ {\it a symmetrized
generalized Cartan matrix}. It obviously defines the matrix
$A$: 
$$ (a_{ij})=\left({2b_{ij}\over b_{ii}}\right)
\tag{1.1.5}
$$ where 
$$ b_{ii} | 2b_{ij}\ \ \text{for any $i,j \in I$}.
\tag{1.1.6} 
$$ 
The matrix $B$ is {\it even primitive}
if there does not exist $n \in \bn$ such that $n>1$ and
the matrix $B/n$ also has all $b_{ij}/n$ integral and all $b_{ii}/n$
even integral.
An even primitive symmetric matrix $(b_{ij})$, $i,j \in I$, with
the properties \thetag{1.1.4} and \thetag{1.1.6} is evidently
equivalent to a symmetrizable generalized Cartan matrix. We call 
primitive $B$ with properties 
\thetag{1.1.4} and \thetag{1.1.6} {\it primitive 
symmetrized generalized Cartan matrix}.

A symmetrizable generalized Cartan matrix $A$ is called {\it hyperbolic}
if its symmetrized generalized Cartan matrix $B$ has exactly one
negative square (equivalently, $A$ has exactly one negative 
eigenvalue). For this case, there exists a {\it geometric realization} of
$A$ (and $B$) in a hyperbolic space. We describe it below.

Let us consider an integral hyperbolic symmetric bilinear form ({\it
lattice}, to be short) which is defined by the symmetric hyperbolic
matrix $B$ modulo kernel. We denote the form by $(\ , \ )$ and its
free $\bz$-module by $M$.  Then the matrix $B$ is the Gram matrix of
the corresponding elements $\alpha_i \in M$, $i \in I$, which generate
the module $M$: $$ B=(b_{ij})=((\alpha_i, \alpha_j)).
\tag{1.1.7}
$$ Then $$ A=(a_{ij})=\left({2(\alpha_i, \alpha_j)\over(\alpha_i,
\alpha_i)}\right).
\tag{1.1.8}
$$ The lattice $M$ is {\it hyperbolic} (i. e., it has exactly one
negative square), {\it even} (i.e. $(x,x)$ is even for any $x \in M$).  
To be shorter, we
will denote a lattice $(M, (\ ,\ ))$ by one letter $M$. The lattice
$M$ defines the cone $$ V(M)=\{x \in M\otimes \br\ |\ (x,x)<0 \}.
\tag{1.1.9}
$$ There exists unique choice of its half $V^+(M)$ such that the cone
$$
\br_+\M =\{x \in V^+(M)\ |\ (x, \alpha_i)\le 0\ \ \text{for any}\  i \in I\} 
\tag{1.1.10} 
$$ is non-empty. Here each inequality $(x, \alpha_i)\le 0$, $i \in I$,
is essential (i. e. the cone will be changed if one removes this
inequality).
 
Using hyperbolic geometry language, the projectivization
$\La^+(M)=V^+(M)/\br_{++}$ is a hyperbolic space, and the
projectivisation $\M$ of the cone $\br_+{\M}$ is a locally finite
polyhedron in the hyperbolic space $\La^+(M)$ with the set of
orthogonal vectors (to faces of $\M$ and directed outwards) $P(\M
)=\{\alpha_i\ |\ i \in I\}$. Since the lattice $M$ is integral
hyperbolic and because of conditions $$ (\alpha_i, \alpha_i) |
2(\alpha_i, \alpha_j) \ \text{for any $i \in I$},
\tag{1.1.11}
$$ $$ (\alpha_i, \alpha_j)\le 0\ \ \text{if $i\not=j \in I$}
\tag{1.1.12} 
$$ (see \thetag{1.1.4} and \thetag{1.1.6}), the polyhedron $\M$ is a
fundamental polyhedron for a discrete reflection group $W$ generated
by reflections $s_{\alpha_i}$ in hyperplanes $$
\Ha_{\alpha_i}=\{ \ \br_{++} x \in \La^+(M)\ |\ 
(x, \alpha_i)=0 \}.
\tag{1.1.13}
$$ It changes places half spaces $\Ha_{\alpha_i}^+$ and
$\Ha_{-\alpha_i}^+$ bounded by the hyperplane $\Ha_{\alpha_i}$ where
$$
\Ha_{\alpha_i}^+=\{ \ \br_{++} x \in \La^+(M)\ |\ 
(x, \alpha_i)\le 0\}.
\tag{1.1.14}
$$ On the hyperbolic lattice $M$ the reflection $s_{\alpha_i}$ is
defined by the formula $$ x \mapsto x-(2(\alpha_i, x)/(\alpha_i,
\alpha_i)) \alpha_i,\ \ x \in M.
\tag{1.1.15}   
$$

Thus, an indecomposable hyperbolic generalized Cartan matrix $A$ is
equivalent to a triplet 
$$ 
(M, W, P(\M)).
\tag{1.1.16} 
$$ 
Here $M$ is a primitive even hyperbolic lattice and 
$W \subset W(M)\subset O^+(M)$ is a reflection subgroup 
generated by
reflections in a set of elements of the lattice $M$ with positive
squares. $W(M)$ denote the group generated by reflections in all 
elements of the lattice $M$ with positive squares, 
$O(M)$ is the orthogonal group of the lattice $M$ and
$O^+(M)$ is its subgroup which fixes the half-cone $V^+(M)$. 
The third object 
$P(\M )\subset M$ in \thetag{1.1.16} 
is an acceptable set of elements orthogonal to the faces
of a fundamental polyhedron $\M\subset \La^+(\M)$ of the group $W$
(the set $P(\M)$ contains exactly one element orthogonal to each face of
$\M$).  The set $P(\M )$ is called {\it acceptable} if it generates
the lattice $M$ and has the property $$ (\alpha, \alpha) |2(\alpha ,
M)\ \ \text{for any $\alpha \in P(\M)$}.
\tag{1.1.17} 
$$ Moreover, we suppose the set $P(\M)$ has an indecomposable Gram
matrix.  Any triplet \thetag{1.1.16} defines a 
hyperbolic generalized Cartan matrix $$ A=\left({2(\alpha,
\alpha^\prime)\over (\alpha, \alpha))}\right), \ \
\alpha , \alpha^\prime \in P(\M).
\tag{1.1.18}
$$ See \cite{N9}, \cite{N10} for details.  We call the triplet
\thetag{1.1.16} a {\it geometric realization} of the hyperbolic
generalized Cartan matrix
\thetag{1.1.18}.  

Suppose that there are natural $\lambda(\alpha )$, $\alpha \in P(\M)$,
such that all $\lambda(\alpha )$ are coprime: 
$$
\text{g.c.d.}(\{\lambda(\alpha)\ |\ \alpha \in P(\M )\})=1, 
\tag{1.1.19} 
$$ 
and elements $\widetilde{\alpha} =\lambda(\alpha )\alpha$ satisfy
the condition
\thetag{1.1.11}, i. e.   
$$ (\widetilde{\alpha} , \widetilde{\alpha}) | 2(\widetilde{\alpha},
\widetilde{\alpha} ^\prime)\ \
\text{for any  
$\alpha, \alpha^\prime \in P(\M)$}.
\tag{1.1.20} 
$$ Equivalently, $$
\lambda(\alpha)(\alpha, \alpha) | 
2\lambda(\alpha^\prime)(\alpha , \alpha^\prime)\ \
\text{for any  
$\alpha, \alpha^\prime \in P(\M)$}.
\tag{1.1.21} 
$$
Then the matrix $$
\widetilde{A}=\left( {2(\widetilde{\alpha},\widetilde{\alpha}^\prime)\over 
(\widetilde{\alpha}, \widetilde{\alpha})}\right) =
\left({2\lambda(\alpha^\prime)(\alpha, \alpha^\prime)\over 
\lambda(\alpha ) (\alpha, \alpha))}\right), \ \ 
\alpha , \alpha^\prime \in P(\M), 
\tag{1.1.22}
$$ is also a hyperbolic generalized Cartan matrix with the same
reflection group $W$. This matrix is called {\it twisted to $A$}.
Coefficients $\lambda (\alpha)$ satisfying \thetag{1.1.19} and
\thetag{1.1.20} (equivalently, \thetag{1.1.19} and 
\thetag{1.1.21}) are called 
{\it twisted coefficients}.  The corresponding generalized Cartan
matrix $\widetilde{A}$ is called {\it twisted to $A$}.  Thus, the
geometric realization of the generalized Cartan matrix $\widetilde{A}$
is equal to $$ (\widetilde{M},\ \widetilde{W},\
\widetilde{P}(\widetilde{\M}))= ([\{\lambda(\alpha)\alpha\ |\ \alpha
\in P(\M )\}]\subset M,\ W, \
\{\lambda(\alpha)\alpha\ |\ \alpha \in P(\M)\}) 
\tag{1.1.23}
$$ where we identify hyperbolic spaces of $M$ and $\widetilde M$. 
(Here $[X]$ denotes a submodule, sublattice generated by the set $X$.) 
Thus, the basic relation between $A$ and its twisted $\widetilde{A}$
is that they have the same reflection (or Weyl) group $W$ and the
same fundamental polyhedron $\M$.

A generalized Cartan matrix is called {\it untwisted} if it is not
twisted to any generalized Cartan matrix different from itself. To
find all possible generalized Cartan matrices, it is sufficient to
find all untwisted generalized Cartan matrices and find for them all
possible sets of twisted coefficients.  For a finite generalized
Cartan matrix $A$ the number of possible twisted coefficients and twisted
to $A$ generalized Cartan matrices $\tilde{A}$ is finite.

Let $A$ be a hyperbolic generalized Cartan matrix and \thetag{1.1.16} 
its geometric realization. We introduce the group $$
\text{Sym\ }(A)=\Sym (P(\M ))=\{g \in O^+(\M )\ |\ g(P(\M ))=P(\M)\ \}
\tag{1.1.24} 
$$ 
which is called {\it group of symmetries} of the generalized Cartan
matrix $A$ (or $P(\M )$). One can consider this group as a subgroup of
$O(M)$ or a subgroup of symmetries of the fundamental polyhedron $\M$.

We use the following definitions from \cite{N10} which are important
for automorphic correction of Lorentzian Kac--Moody algebras
corresponding to hyperbolic generalized Cartan matrices. One has a
chance to find an automorphic correction only for cases described
below. See Sects. 2.4, 2.5. 

\definition{Definition 1.1.1} A hyperbolic generalized Cartan matrix $A$ 
has {\it restricted arithmetic type} if it is not empty and the
semi-direct product 
$$ W \rtimes \Sym (A)=W \rtimes \Sym (P(\M ))
\tag{1.1.25}
$$ has a finite index in $O^+(M)$.
\enddefinition

\definition{Definition 1.1.2} A hyperbolic generalized Cartan matrix $A$ 
has a {\it lattice Weyl vector} if there exists an element $\rho \in
M\otimes \bq$ such that for a constant $N>0$ 
$$ (\rho,
\alpha)=-(\alpha, \alpha )/2\ \ \text{for any $\alpha \in P(\M )$}.
\tag{1.1.26}
$$ Then $\rho $ is called a {\it lattice Weyl vector}.

More generally, $A$ has a {\it generalized lattice Weyl vector} if
there exists a non-zero element $\rho \in M \otimes \bq$ such that for
a constant $N>0$ one has $$ 0\le -(\rho, \alpha)\le N \ \ \text{for
any $\alpha \in P(\M )$}.
\tag{1.1.27} 
$$ The element $\rho$ is called a {\it generalized lattice Weyl
vector}.
\enddefinition

One can prove (see \cite{N10}) that if $\rk A=\rk M \ge 3$
(equivalently, $\dim \La^+(M)\ge 2$) and $A$ has a restricted
arithmetic type, then $\br_{++}\rho \in \M$ and $(\rho, \rho )\le 0$.
Therefore, one should consider two cases:

\definition{Definition 1.1.3} A hyperbolic generalized Cartan matrix 
has {\it elliptic type} if it has a restricted arithmetic type and has
a generalized lattice Weyl vector $\rho $ with negative square:
$(\rho, \rho )<0$. This is equivalent to finiteness of index
$[O(M):W]< \infty $ or finiteness of volume of $\M$ (in particular,
the set $P(\M)$ is finite).

A hyperbolic generalized Cartan matrix $A$ has {\it parabolic type} if
it has restricted arithmetic type and has a generalized lattice Weyl
vector $\rho $ with zero square: $(\rho, \rho )=0$, and it 
does not have a generalized lattice Weyl vector with 
negative square (i.e. it is
not elliptic).  For this case the group of symmetries $\Sym
(A)=\Sym(P(\M))\subset O(\rho^\perp_M)$ is a crystallographic group
and is Abelian up to a finite index.
\enddefinition

It is proved in \cite{N10} (in fact, it was done in \cite{N3},
\cite{N4}, \cite{N5}) that number of these cases is essentially finite. 

The case when $A$ has elliptic or parabolic type and at the same time
has a lattice Weyl vector is especially important from the point of
view of the theory of Lorentzian Kac--Moody algebras. 
See Sects 2.4, 2.5. It is proved in
\cite{N10} that number of hyperbolic generalized Cartan matrices of
elliptic type and with a lattice Weyl vector is finite for rank $\ge
3$.

In the next Section we want to classify all symmetric generalized
Cartan matrices of elliptic type and of rank $3$ 
with a lattice Weyl vector. 
More generally, we will describe all symmetrized
generalized Cartan matrices of elliptic type with a lattice Weyl
vector which are twisted to symmetric generalized Cartan matrices of
rank 3 (the last may does not have a lattice Weyl vector).

\subhead
1.2. The classification of generalized Cartan matrices of elliptic type 
and of rank $3$ with a lattice Weyl vector which are twisted to symmetric
generalized Cartan matrices
\endsubhead

Let $A$ be a generalized Cartan matrix of elliptic type.  
Suppose that $A$ is
twisted to a symmetric generalized Cartan matrix $\overline{A}$.
Obviously, then $\overline{A}$ is also elliptic.  Let us take a
geometric realization $G(A)=(M, W, P(\M))$ of $A$.  Then for any
$\alpha \in P(\M)$ one has $\alpha =\lambda (\alpha )\delta (\alpha)$
where $(\delta (\alpha), \delta (\alpha))=2$, $\lambda (\alpha ) \in
\bn$ and the matrix $\overline{A}=(\delta (\alpha ), \delta
(\alpha^\prime)),\ \
\alpha, \alpha^\prime \in P(\M)$,  
is a symmetric generalized Cartan matrix.  Suppose additionally that
$\rk A=3$. Then $\M$ is a polygon $A_1A_2....A_n$ on a hyperbolic
plane, and we can numerate elements of $P(\M)$ and
$\overline{P}(\M)=\{ \delta (\alpha)=\alpha/\lambda (\alpha)\ |\
\alpha \in P(\M )\}$ 
by $1,...,n$ where $\alpha_i$ is orthogonal to the side $A_iA_{i+1}$
of the polygon. We naturally denote $\delta_i=\delta (\alpha_i)$ and
$\lambda_i=\lambda (\alpha_i)$. Then $A$ and its geometric realization
are equivalent to $(1+[n/2])\times n$ matrix $G(A)$. It has the first
line $$
\lambda_1,...,\lambda_{n}
$$ and it has $i+1$-th line $$ -(\delta_1,\delta_{1+i}), -(\delta_2,
\delta_{2+i}),...,-(\delta_n, \delta_{n+i}), $$ $1\le i \le [n/2]$,
where we consider the second index $j$, $1 \le j \le n$, modulo $n$.
Thus, a $j$-th column of $G(A)$ is 
$$ (\lambda_j, (\delta_j,
\delta_{j+1}), (\delta_j, \delta_{j+2}),..., (\delta_j,
\delta_{j+[n/2]}))^t.
$$ 
We want to find all matrices $G(A)$ of this
type having a lattice Weyl vector $\rho$ which is an element $\rho
\in M\otimes \bq$ satisfying \thetag{1.1.26}. Equivalently, 
$$ 
(\rho, \delta_i)=-\lambda_i,\ \ i=1,...,n.
\tag{1.2.1}
$$ 
The geometric realization $G(A)$ is convenient because it easily
permits to find the group of symmetries $\Sym (A)$, it has smaller
coefficients and smaller size than the generalized Cartan matrix $A$
and the corresponding symmetric generalized Cartan matrix $B$.

We want to prove the following classification

\proclaim{Theorem 1.2.1} All geometric realizations $G(A)$ of 
hyperbolic generalized Cartan matrices $A$ of elliptic type with a
lattice Weyl vector and with properties: $A$ has rank $3$, the matrix
$A$ is twisted to a symmetric generalized Cartan matrix, all twisting
coefficients $\lambda_i$ satisfy the inequality $\lambda_i \le 12$,
are given in Table 1 below.
\endproclaim

\proclaim{Conjecture 1.2.2} Table 1 gives the complete 
list of hyperbolic generalized Cartan matrices $A$ of elliptic type
with a lattice Weyl vector and with the properties: $A$ has rank $3$ and
the matrix $A$ is twisted to a symmetric generalized Cartan matrix 
(i.e. one can drop the inequalities $\lambda_i \le 12$ from the Theorem
1.2.1).
\endproclaim

We have the following arguments to support Conjecture 1.2.2. First, it
is known \cite{N10} that the number of all hyperbolic generalized Cartan
matrices of elliptic type with a lattice Weyl vector is finite for
rank $\ge 3$. Thus, there exists an absolute constant $m$ such that
all $\lambda_i \le m$.  Second, we did our calculations for all
$\lambda_i \le 12$, but the result of these calculations has only
matrices with all $\lambda_i \le 6$. Thus, in between $6$ and $12$
there are no new solutions.  In further publications we hope to
present additional arguments which prove Conjecture 1.2.2. 

\smallpagebreak 

\centerline{\bf Table 1.} 

\smallpagebreak

\centerline{Geometric realizations of hyperbolic generalized Cartan matrices}
\centerline{of elliptic type 
with a lattice Weyl vector which are twisted}  
\centerline{to symmetric generalized Cartan matrices and have rank 3.}  

\smallpagebreak

\noindent
$r=-59/2: \ \ G(A)=\left(\matrix {1}&{2}&{2}\cr {}&{}&{}\cr
{0}&{1}&{2}\cr
\endmatrix\right). $ 
\ \ \ \ 
$r=-22: \ \ G(A)=\left(\matrix {2}&{1}&{1}\cr {}&{}&{}\cr
{0}&{1}&{2}\cr
\endmatrix\right) .
$

\noindent
$r=-16:\ \ G(A)=\left(\matrix {1}&{4}&{2}\cr {}&{}&{}\cr
{0}&{2}&{2}\cr
\endmatrix\right) .$
\ \ \ \ 
$r=-23/2:\ \ G(A)=\left(\matrix {1}&{1}&{1}\cr {}&{}&{}\cr
{0}&{1}&{2}\cr
\endmatrix\right) .$

\noindent  
$r=-10:\ \ G(A)=\left(\matrix {1}&{2}&{2}\cr {}&{}&{}\cr
{0}&{2}&{2}\cr
\endmatrix\right) .$
\ \ \ \ 
$r=-17/2:\ \ G(A)=\left(\matrix {2}&{2}&{1}\cr {}&{}&{}\cr
{0}&{2}&{2}\cr
\endmatrix\right) .$

\noindent  
$r=-7:\ \ G(A)=
\left(\matrix
{1}&{1}&{2}\cr {}&{}&{}\cr {0}&{2}&{2}\cr
\endmatrix\right),
\ \
\left(\matrix
{1}&{3}&{3}&{1}\cr {}&{}&{}&{}\cr {0}&{1}&{0}&{1}\cr
{3}&{3}&{3}&{3}\cr
\endmatrix\right) .
$
\ \ \ 
$r=-6:\ \ G(A)=
\left(\matrix
{1}&{6}&{3}&{2}\cr {}&{}&{}&{}\cr {0}&{2}&{0}&{2}\cr
{3}&{6}&{3}&{6}\cr
\endmatrix\right) . 
$

\noindent
$r=-11/2:\ \ G(A)=\left(\matrix {1}&{2}&{1}\cr {}&{}&{}\cr
{0}&{2}&{2}\cr
\endmatrix\right),
\ \
\left(\matrix
{2}&{2}&{1}\cr {}&{}&{}\cr {1}&{2}&{2}\cr
\endmatrix\right) .
$
\ \ \ \ 
$r=-4:\ \ G(A)=\left(\matrix {1}&{1}&{1}\cr {}&{}&{}\cr {1}&{1}&{2}\cr
\endmatrix\right),
$

\noindent
$
\left(\matrix
{1}&{1}&{2}\cr {}&{}&{}\cr {1}&{2}&{2}\cr
\endmatrix\right),
\ \ 
\left(\matrix
{1}&{2}&{2}\cr {}&{}&{}\cr {2}&{2}&{2}\cr
\endmatrix\right),
\ \ 
\left(\matrix
{1}&{4}&{4}&{1}\cr {}&{}&{}&{}\cr {0}&{2}&{0}&{2}\cr
{4}&{4}&{4}&{4}\cr
\endmatrix\right),
\ \ 
\left(\matrix
{1}&{4}&{4}&{2}\cr {}&{}&{}&{}\cr {0}&{2}&{2}&{2}\cr
{4}&{6}&{4}&{6}\cr
\endmatrix\right). 
$

\noindent
$r=-7/2:\ \ G(A)=
\left(\matrix
{1}&{1}&{1}\cr {}&{}&{}\cr {0}&{2}&{2}\cr
\endmatrix\right) .
$
\ \ \ \ 
$r=-5/2: \ \ G(A)=
\left(\matrix
{1}&{2}&{1}\cr {}&{}&{}\cr {2}&{2}&{2}\cr
\endmatrix\right) .
$

\noindent
$r=-13/6:\ \ G(A)=
\left(\matrix
{1}&{1}&{1}\cr {}&{}&{}\cr {1}&{2}&{2}\cr
\endmatrix\right) .
$
\ \ \ \ 
$r=-17/8:\ \ G(A)=
\left(\matrix
{1}&{1}&{2}&{2}\cr {}&{}&{}&{}\cr {0}&{2}&{0}&{2}\cr
{4}&{4}&{4}&{4}\cr
\endmatrix\right) .
$

\noindent
$r=-2:\ G(A)=
\left(\matrix
{1}&{1}&{1}&{1}\cr
{}&{}&{}&{}\cr
{0}&{1}&{0}&{1}\cr
{3}&{3}&{3}&{3}\cr
\endmatrix\right), 
\  
\left(\matrix
{1}&{2}&{1}&{2}\cr
{}&{}&{}&{}\cr
{0}&{2}&{0}&{2}\cr
{3}&{6}&{3}&{6}\cr
\endmatrix\right), 
\  
\left(\matrix
{1}&{2}&{2}&{1}\cr
{}&{}&{}&{}\cr
{0}&{2}&{0}&{2}\cr
{4}&{4}&{4}&{4}\cr
\endmatrix\right) ,
\   
\left(\matrix
{1}&{2}&{2}&{2}\cr
{}&{}&{}&{}\cr
{0}&{2}&{2}&{2}\cr
{4}&{6}&{4}&{6}\cr
\endmatrix\right) . 
$

\noindent
$r=-3/2:\ \   
G(A)=
\left(\matrix
{1}&{1}&{1}\cr
{}&{}&{}\cr
{2}&{2}&{2}\cr
\endmatrix\right),
\ \ 
\left(\matrix
{2}&{3}&{3}&{2}&{1}\cr
{}&{}&{}&{}&{}\cr
{0}&{2}&{0}&{2}&{2}\cr
{6}&{6}&{6}&{7}&{6}\cr
\endmatrix\right).
$

\noindent
$r=-1:\ \ G(A)= 
\left(\matrix
{1}&{1}&{1}&{1}\cr
{}&{}&{}&{}\cr
{0}&{2}&{0}&{2}\cr
{4}&{4}&{4}&{4}\cr
\endmatrix\right),
\ \  
\left(\matrix
{1}&{1}&{1}&{2}\cr
{}&{}&{}&{}\cr
{0}&{2}&{2}&{2}\cr
{4}&{6}&{4}&{6}\cr
\endmatrix\right), 
\ \  
\left(\matrix
{1}&{1}&{1}&{1}\cr
{}&{}&{}&{}\cr
{1}&{1}&{1}&{1}\cr
{4}&{4}&{4}&{4}\cr
\endmatrix\right),  
$

\noindent
$  
\left(\matrix
{1}&{2}&{1}&{2}\cr
{}&{}&{}&{}\cr
{2}&{2}&{2}&{2}\cr
{4}&{10}&{4}&{10}\cr
\endmatrix\right), 
\ \ 
\left(\matrix
{1}&{1}&{2}&{2}\cr
{}&{}&{}&{}\cr
{2}&{2}&{2}&{2}\cr
{6}&{6}&{6}&{6}\cr
\endmatrix\right), 
\ \ 
\left(\matrix
{1}&{3}&{4}&{3}&{1}&{3}&{4}&{3}\cr
{}&{}&{}&{}&{}&{}&{}&{}\cr
{0}&{0}&{0}&{0}&{0}&{0}&{0}&{0}\cr
{4}&{4}&{4}&{14}&{4}&{4}&{4}&{14}\cr
{6}&{6}&{24}&{24}&{6}&{6}&{24}&{24}\cr
{4}&{20}&{34}&{20}&{4}&{20}&{34}&{20}\cr
\endmatrix\right) .
$

\noindent
$r=-2/3:\ \  
G(A)=
\left(\matrix
{1}&{1}&{1}&{1}\cr
{}&{}&{}&{}\cr
{1}&{2}&{1}&{2}\cr
{5}&{5}&{5}&{5}\cr
\endmatrix\right) .
$
\ \ 
$
r=-5/8:\ \  
G(A)=\left(\matrix
{1}&{2}&{1}&{2}&{1}\cr
{}&{}&{}&{}&{}\cr
{0}&{2}&{2}&{0}&{2}\cr
{4}&{14}&{4}&{8}&{8}\cr
\endmatrix\right) .
$

\noindent
$r=-1/2:\ \ G(A)=
\left(\matrix
{1}&{1}&{1}&{1}\cr
{}&{}&{}&{}\cr
{2}&{2}&{2}&{2}\cr
{6}&{6}&{6}&{6}\cr
\endmatrix\right),
\ \   
\left(\matrix
{2}&{2}&{1}&{1}&{1}\cr
{}&{}&{}&{}&{}\cr
{0}&{2}&{2}&{2}&{2}\cr
{6}&{10}&{6}&{10}&{6}\cr
\endmatrix\right), 
\ \  
\left(\matrix
{1}&{1}&{1}&{1}&{1}&{1}\cr
{}&{}&{}&{}&{}&{}\cr
{0}&{0}&{0}&{0}&{0}&{0}\cr
{4}&{4}&{4}&{4}&{4}&{4}\cr
{6}&{6}&{6}&{6}&{6}&{6}\cr
\endmatrix\right), 
$

\noindent
$     
\left(\matrix
{2}&{2}&{1}&{1}&{2}&{2}\cr
{}&{}&{}&{}&{}&{}\cr
{0}&{2}&{2}&{2}&{0}&{2}\cr
{6}&{10}&{10}&{6}&{8}&{8}\cr
{10}&{14}&{10}&{10}&{14}&{10}\cr
\endmatrix\right), 
\ \ \   
\left(\matrix
{2}&{2}&{1}&{2}&{2}&{1}\cr
{}&{}&{}&{}&{}&{}\cr
{0}&{2}&{2}&{0}&{2}&{2}\cr
{6}&{16}&{6}&{6}&{16}&{6}\cr
{18}&{18}&{6}&{18}&{18}&{6}\cr
\endmatrix\right), 
$

\noindent
$  
\left(\matrix
{2}&{2}&{1}&{2}&{2}&{2}&{2}\cr
{}&{}&{}&{}&{}&{}&{}\cr
{0}&{2}&{2}&{0}&{2}&{0}&{2}\cr
{6}&{16}&{6}&{8}&{8}&{8}&{8}\cr
{18}&{18}&{10}&{14}&{16}&{14}&{10}\cr
\endmatrix\right) . 
$
\ \ 
$r=-2/5:\ \   
G(A)= 
\left(\matrix
{1}&{2}&{2}&{1}&{2}&{2}\cr
{}&{}&{}&{}&{}&{}\cr
{2}&{2}&{2}&{2}&{2}&{2}\cr
{8}&{8}&{18}&{8}&{8}&{18}\cr
{7}&{22}&{22}&{7}&{22}&{22}\cr
\endmatrix\right) .$

\noindent    
$r=-7/18:\ \    
G(A)=
\left(\matrix
{1}&{1}&{1}&{1}&{1}\cr
{}&{}&{}&{}&{}\cr
{0}&{2}&{0}&{2}&{2}\cr
{6}&{6}&{6}&{7}&{6}\cr
\endmatrix\right) .
$

\noindent  
$r=-1/4:\ \  
G(A)=
\left(\matrix
{1}&{1}&{1}&{1}&{1}&{1}\cr
{}&{}&{}&{}&{}&{}\cr
{1}&{1}&{1}&{1}&{1}&{1}\cr
{7}&{7}&{7}&{7}&{7}&{7}\cr
{10}&{10}&{10}&{10}&{10}&{10}\cr
\endmatrix\right), 
\ \  
\left(\matrix
{1}&{2}&{1}&{1}&{2}&{1}\cr
{}&{}&{}&{}&{}&{}\cr
{2}&{2}&{2}&{2}&{2}&{2}\cr
{6}&{14}&{14}&{6}&{14}&{14}\cr
{10}&{34}&{10}&{10}&{34}&{10}\cr
\endmatrix\right) .
$

\noindent 
$r=-2/9:
\ \  
G(A)=
\left(\matrix
{1}&{2}&{2}&{1}&{1}&{2}&{2}&{1}\cr
{}&{}&{}&{}&{}&{}&{}&{}\cr
{0}&{2}&{0}&{2}&{0}&{2}&{0}&{2}\cr
{6}&{6}&{12}&{12}&{6}&{6}&{12}&{12}\cr
{7}&{18}&{34}&{18}&{7}&{18}&{34}&{18}\cr
{11}&{38}&{38}&{11}&{11}&{38}&{38}&{11}\cr
\endmatrix\right) .
$

\noindent
$r=-1/6:\ \  
G(A)=
\left(\matrix
{1}&{1}&{1}&{1}&{1}&{1}\cr
{}&{}&{}&{}&{}&{}\cr
{2}&{2}&{2}&{2}&{2}&{2}\cr
{10}&{10}&{10}&{10}&{10}&{10}\cr
{14}&{14}&{14}&{14}&{14}&{14}\cr
\endmatrix\right),
\ \ 
\left(\matrix
{2}&{2}&{1}&{1}&{1}&{1}&{1}\cr
{}&{}&{}&{}&{}&{}&{}\cr
{1}&{2}&{2}&{2}&{2}&{2}&{2}\cr
{8}&{16}&{10}&{10}&{10}&{16}&{8}\cr
{22}&{26}&{14}&{14}&{26}&{22}&{10}\cr
\endmatrix\right),
$ 

\noindent
$
\left(\matrix
{2}&{2}&{1}&{1}&{1}&{1}&{2}&{2}\cr
{}&{}&{}&{}&{}&{}&{}&{}\cr
{1}&{2}&{2}&{2}&{2}&{2}&{1}&{2}\cr
{8}&{16}&{10}&{10}&{16}&{8}&{11}&{11}\cr
{22}&{26}&{14}&{26}&{22}&{16}&{23}&{16}\cr
{26}&{22}&{22}&{26}&{26}&{22}&{22}&{26}\cr
\endmatrix\right), 
\ \ 
\left(\matrix
{2}&{2}&{1}&{1}&{1}&{2}&{2}&{1}\cr
{}&{}&{}&{}&{}&{}&{}&{}\cr
{1}&{2}&{2}&{2}&{2}&{1}&{2}&{2}\cr
{8}&{16}&{10}&{16}&{8}&{8}&{25}&{8}\cr
{22}&{26}&{26}&{22}&{10}&{37}&{37}&{10}\cr
{26}&{46}&{26}&{14}&{26}&{46}&{26}&{14}\cr
\endmatrix\right),
$

\noindent
$ 
\left(\matrix
{2}&{2}&{1}&{1}&{2}&{2}&{1}&{1}\cr
{}&{}&{}&{}&{}&{}&{}&{}\cr
{1}&{2}&{2}&{2}&{1}&{2}&{2}&{2}\cr
{8}&{16}&{16}&{8}&{8}&{16}&{16}&{8}\cr
{22}&{47}&{22}&{10}&{22}&{47}&{22}&{10}\cr
{50}&{50}&{14}&{14}&{50}&{50}&{14}&{14}\cr
\endmatrix\right), 
\ \ 
\left(\matrix
{2}&{2}&{1}&{1}&{1}&{2}&{2}&{2}&{2}\cr
{}&{}&{}&{}&{}&{}&{}&{}&{}\cr
{1}&{2}&{2}&{2}&{2}&{1}&{2}&{1}&{2}\cr
{8}&{16}&{10}&{16}&{8}&{11}&{11}&{11}&{11}\cr
{22}&{26}&{26}&{22}&{16}&{23}&{25}&{23}&{16}\cr
{26}&{46}&{26}&{26}&{22}&{37}&{37}&{22}&{26}\cr
\endmatrix\right), 
$ 

\noindent
$
\left(\matrix
{2}&{2}&{1}&{1}&{2}&{2}&{1}&{2}&{2}\cr
{}&{}&{}&{}&{}&{}&{}&{}&{}\cr
{1}&{2}&{2}&{2}&{1}&{2}&{2}&{1}&{2}\cr
{8}&{16}&{16}&{8}&{8}&{25}&{8}&{11}&{11}\cr
{22}&{47}&{22}&{10}&{37}&{37}&{16}&{23}&{16}\cr
{50}&{50}&{14}&{26}&{46}&{47}&{22}&{22}&{26}\cr
\endmatrix\right),
$

\noindent
$ 
\left(\matrix
{2}&{2}&{1}&{2}&{2}&{1}&{2}&{2}&{1}\cr
{}&{}&{}&{}&{}&{}&{}&{}&{}\cr
{1}&{2}&{2}&{1}&{2}&{2}&{1}&{2}&{2}\cr
{8}&{25}&{8}&{8}&{25}&{8}&{8}&{25}&{8}\cr
{37}&{37}&{10}&{37}&{37}&{10}&{37}&{37}&{10}\cr
{46}&{26}&{26}&{46}&{26}&{26}&{46}&{26}&{26}\cr
\endmatrix\right), 
$

\noindent
$
\left(\matrix
{2}&{2}&{1}&{1}&{2}&{2}&{2}&{2}&{2}&{2}\cr
{}&{}&{}&{}&{}&{}&{}&{}&{}&{}\cr
{1}&{2}&{2}&{2}&{1}&{2}&{1}&{2}&{1}&{2}\cr
{8}&{16}&{16}&{8}&{11}&{11}&{11}&{11}&{11}&{11}\cr
{22}&{47}&{22}&{16}&{23}&{25}&{23}&{25}&{23}&{16}\cr
{50}&{50}&{26}&{22}&{37}&{37}&{37}&{37}&{22}&{26}\cr
{47}&{46}&{26}&{26}&{46}&{47}&{46}&{26}&{26}&{46}\cr
\endmatrix\right), 
$  

\noindent
$
\left(\matrix
{2}&{2}&{1}&{2}&{2}&{1}&{2}&{2}&{2}&{2}\cr
{}&{}&{}&{}&{}&{}&{}&{}&{}&{}\cr
{1}&{2}&{2}&{1}&{2}&{2}&{1}&{2}&{1}&{2}\cr
{8}&{25}&{8}&{8}&{25}&{8}&{11}&{11}&{11}&{11}\cr
{37}&{37}&{10}&{37}&{37}&{16}&{23}&{25}&{23}&{16}\cr
{46}&{26}&{26}&{46}&{47}&{22}&{37}&{37}&{22}&{47}\cr
{26}&{46}&{26}&{50}&{50}&{26}&{46}&{26}&{50}&{50}\cr
\endmatrix\right), 
$

\noindent
$ 
\left(\matrix
{2}&{2}&{1}&{2}&{2}&{2}&{2}&{1}&{2}&{2}\cr
{}&{}&{}&{}&{}&{}&{}&{}&{}&{}\cr
{1}&{2}&{2}&{1}&{2}&{1}&{2}&{2}&{1}&{2}\cr
{8}&{25}&{8}&{11}&{11}&{8}&{25}&{8}&{11}&{11}\cr
{37}&{37}&{16}&{23}&{16}&{37}&{37}&{16}&{23}&{16}\cr
{46}&{47}&{22}&{22}&{47}&{46}&{47}&{22}&{22}&{47}\cr
{50}&{50}&{14}&{50}&{50}&{50}&{50}&{14}&{50}&{50}\cr
\endmatrix\right),
$  

\noindent
$
\left(\matrix
{2}&{2}&{1}&{2}&{2}&{2}&{2}&{2}&{2}&{2}&{2}\cr
{}&{}&{}&{}&{}&{}&{}&{}&{}&{}&{}\cr
{1}&{2}&{2}&{1}&{2}&{1}&{2}&{1}&{2}&{1}&{2}\cr
{8}&{25}&{8}&{11}&{11}&{11}&{11}&{11}&{11}&{11}&{11}\cr
{37}&{37}&{16}&{23}&{25}&{23}&{25}&{23}&{25}&{23}&{16}\cr
{46}&{47}&{22}&{37}&{37}&{37}&{37}&{37}&{37}&{22}&{47}\cr
{50}&{50}&{26}&{46}&{47}&{46}&{47}&{46}&{26}&{50}&{50}\cr
\endmatrix\right), 
$

\noindent
$ 
\left(\matrix
{1}&{2}&{2}&{1}&{2}&{2}&{1}&{2}&{2}&{1}&{2}&{2}\cr
{}&{}&{}&{}&{}&{}&{}&{}&{}&{}&{}&{}\cr
{0}&{0}&{0}&{0}&{0}&{0}&{0}&{0}&{0}&{0}&{0}&{0}\cr
{4}&{4}&{14}&{4}&{4}&{14}&{4}&{4}&{14}&{4}&{4}&{14}\cr
{6}&{24}&{24}&{6}&{24}&{24}&{6}&{24}&{24}&{6}&{24}&{24}\cr
{20}&{34}&{20}&{20}&{34}&{20}&{20}&{34}&{20}&{20}&{34}&{20}\cr
{24}&{24}&{48}&{24}&{24}&{48}&{24}&{24}&{48}&{24}&{24}&{48}\cr
{14}&{50}&{50}&{14}&{50}&{50}&{14}&{50}&{50}&{14}&{50}&{50}\cr
\endmatrix\right) .
$

\noindent
$r=-1/8:\ \
G(A)=\left(\matrix
{1}&{1}&{1}&{1}&{1}&{1}&{1}&{1}\cr
{}&{}&{}&{}&{}&{}&{}&{}\cr
{0}&{2}&{0}&{2}&{0}&{2}&{0}&{2}\cr
{8}&{8}&{8}&{8}&{8}&{8}&{8}&{8}\cr
{14}&{16}&{14}&{16}&{14}&{16}&{14}&{16}\cr
{18}&{18}&{18}&{18}&{18}&{18}&{18}&{18}\cr
\endmatrix\right) .
$

\noindent
$r=-1/24:\ \
G(A)=
\left(\matrix
{1}&{1}&{1}&{1}&{1}&{1}&{1}&{1}&{1}&{1}&{1}&{1}\cr
{}&{}&{}&{}&{}&{}&{}&{}&{}&{}&{}&{}\cr
{1}&{2}&{1}&{2}&{1}&{2}&{1}&{2}&{1}&{2}&{1}&{2}\cr
{11}&{11}&{11}&{11}&{11}&{11}&{11}&{11}&{11}&{11}&{11}&{11}\cr
{23}&{25}&{23}&{25}&{23}&{25}&{23}&{25}&{23}&{25}&{23}&{25}\cr
{37}&{37}&{37}&{37}&{37}&{37}&{37}&{37}&{37}&{37}&{37}&{37}\cr
{46}&{47}&{46}&{47}&{46}&{47}&{46}&{47}&{46}&{47}&{46}&{47}\cr
{50}&{50}&{50}&{50}&{50}&{50}&{50}&{50}&{50}&{50}&{50}&{50}\cr
\endmatrix\right) .
$

\demo{Proof of Theorem 1.2.1} We fix $m\in \bn$. 
We give a finite algorithm which 
permits to find all matrices $G(A)$ of our type with  
$\lambda_i \le m$ for all $1\le \lambda_i \le n$. 

{\bf Step 1.}  
It is known \cite{N4}, that 
there exist three consecutive sides of the 
polygon $A_1...A_n$ which we denote by $A_1A_2$, $A_2A_3$, $A_3A_4$ such 
that for the orthogonal vectors $\delta_1, \delta_2, \delta_3 \in 
\overline{P}(\M)$ to 
these sides one has 
$$
-2 \le (\delta_1,\delta_2) \le 0;\ \ 
-14< (\delta_1,\delta_3)\le 0; \ \ 
-2 \le (\delta_2, \delta_3) \le 0.
\tag{1.2.2}
$$
We remark that since 
$$
(\delta_i, \delta_i)=2
\tag{1.2.3}
$$ 
for any $i$, one should have 
$$
-2\le (\delta_i, \delta_{i+1})\le 0
\tag{1.2.4}
$$ 
since the lines $A_iA_{i+1}$ and 
$A_{i+1}A_{i+2}$ should have a common point and the polygon $A_1...A_n$ 
has acute ($\le \pi /2$) angles. The last is equivalent to 
$$
(\delta_i, \delta_j)\le 0\ \ \text{if $i\not=j$}.
\tag{1.2.5}
$$  
Obviously, one can find all (in a finite set)  
possibilities for the symmetric hyperbolic integral matrices  
$g=(\delta_i, \delta_j)$, $1\le i,j \le 3$, satisfying \thetag{1.2.2}.  

We fix one of these $g$. It defines 
the rational hyperbolic bilinear form $M\otimes \bq$. 
For the fixed $g$ 
we find all natural $(\lambda_1, \lambda_2, \lambda_3)$ 
such that 
$$
\lambda_i \le m
\tag{1.2.6} 
$$
and $\lambda_i^2(\delta_i, \delta_i)\ |\ 
2\lambda_i\lambda_j (\delta_i, \delta_j)$, equivalently  
$$
\lambda_i(\delta_i, \delta_i)\ |\ 2\lambda_j (\delta_i, \delta_j)
\tag{1.2.7}
$$
for all $1 \le i,j \le 3$. Further we consider only these 
$(g,\ (\lambda_1,\lambda_2,\lambda_3))$.  
For the fixed 
$(g,\ (\lambda_1, \lambda_2, \lambda_3))$, there exists a unique 
$\rho \in M\otimes \bq$ (a lattice Weyl vector) 
such that $(\rho, \delta_i)=-\lambda_i$, and we calculate 
$$
r=(\rho, \rho).
\tag{1.2.8}
$$ 
Considering all $(g,\ (\lambda_1, \lambda_2, \lambda_3))$ above and 
calculating $r$ for them,  
we find a finite set $R$ of all such $r$ with $r < 0$.   

{\bf Step 2.}  
We fix $r \in R$. Like above, we find all possible hyperbolic 
Gram matrices 
$g_i$ of 
$\delta_i$, $\delta_{i+1}$,\ $\delta_{i+2}$ and twisting 
coefficients $\Lambda_i=(\lambda_i, \lambda_{i+1}, \lambda_{i+2})$ 
which satisfy the conditions \thetag{1.2.3}---\thetag{1.2.7} and  
calculations as  
above of the lattice Weyl vector $\rho$ give the fixed square 
$(\rho, \rho)=r$ (i.e. \thetag{1.2.8} is valid for the fixed $r$ where 
we calculate $\rho$ using $\delta_i,\ \delta_{i+1}, \delta_{i+2}$). 
For fixed 
$\lambda_i$, $\lambda_{i+1}$, $\lambda_{i+2}$ the set of possible such 
Gram matrices $g_i$ is finite. Really, geometrically it is clear that if 
$\lambda_i$, $\lambda_{i+1}$, $\lambda_{i+2}$ and 
$(\delta_i, \delta_{i+1})$, $(\delta_{i+1}, \delta_{i+2})$ are fixed, 
then the function $(\rho, \rho)$ will be increasing for decreasing 
(negative) $(\delta_i, \delta_{i+2})$. Thus, we get a finite set $K^3_r$ 
of possibilities for $(g_i, \Lambda_i)$. 
Actually, one can prove that the set $K^{3}_r$ will be also finite  
if we drop inequalities $\lambda_i\le m$. 
  
Consider datum $(g, \Lambda)$ from $K^3_r$ such that $-2\le g_{i(i+2)}$. 
The ones with coprime $\Lambda$ give all $G(A)$ with $n=3$ we are 
looking for (with the fixed invariant $r$). 
We remove all datum with $-2\le g_{i(i+2)}$ from $K^3_r$ 
and denote the rest by $(K^3_r)^\prime$.    

{\bf Step 3.}  
Using the described set $(K^3_r)^\prime$, we find all Gram matrices 
$$
(g_{ij})=(\delta_i, \delta_j),\ \ 1\le i,j \le 4,
\tag{1.2.9}
$$
and $\Lambda=(\lambda_1, \lambda_2, \lambda_3, \lambda_4)$ 
such that $(g_{ij})$ is hyperbolic  
of rank $3$ and the properties  \thetag{1.2.4} for $1 \le i\le 3$, 
\thetag{1.2.3}, \thetag{1.2.5}, \thetag{1.2.6}, \thetag{1.2.7} 
for $1\le i, j \le 4$, and \thetag{1.2.8} for main submatrices 
$$
(g_{ij})=(\delta_i, \delta_j), \ \ 1\le i,j \le 3, 
\tag{1.2.10}
$$
and 
$$
(g_{ij})=(\delta_i, \delta_j), \ \ 2\le i,j \le 4, 
\tag{1.2.11}
$$
are valid. Equivalently, \thetag{1.2.10} with 
$(\lambda_1, \lambda_2, \lambda_3)$ and \thetag{1.2.11} with 
$(\lambda_2, \lambda_3, \lambda_4)$ both belong to $(K^3_r)^\prime$. 

The last condition defines finite set of possibilities for 
all coefficients of the matrix \thetag{1.2.9} except the 
coefficients $g_{14}=g_{41}=(\delta_1, \delta_4)$. Using \thetag{1.2.8} 
for the matrix $(g_{ij})$, $1\le i,j \le 3$, we can calculate 
the coefficient $g_{14}$ using other coefficients. Really, elements 
$\delta_1$, $\delta_2$, $\rho$ generate over $\bq$ the 
same 3-dimensional lattice as $\delta_1$, $\delta_2$, $\delta_3$. Thus, 
if we know $g_{24}=(\delta_2, \delta_4)$, $g_{34}=(\delta_3, \delta_4)$,  
$(\rho, \delta_4)=-\lambda_4$, we can find 
$g_{14}=(\delta_1,\delta_4)$.  
Further, we should check that $g_{14}$ is integral,  
$g_{14}\le 0$, the condition \thetag{1.2.7} is satisfied 
for $g_{14}=(\delta_1,\delta_4)$,   
and the determinant of the matrix 
\thetag{1.2.9} is equal to $0$. Thus, we can find 
the finite set of matrices \thetag{1.2.9} and $\Lambda$.  
We denote this set by $K^4_r$. 

Consider datum $(g, \Lambda)$ from $K^4_r$ such that $-2\le g_{14}$. 
Those with coprime $\Lambda$ give all $G(A)$ with $n=4$ we are 
looking for. We remove all datum with $-2\le g_{14}$ from $K^4_r$ 
and denote the rest by $(K^4_r)^\prime$. 

{\bf Step $n-1$.}  
Now we define the sets $K^n_r$ and $(K^n_r)^\prime$, 
$n \ge 5$. The set $K^n_r$ is the set of all pairs 
$(g, \Lambda)$ where $g$ is $n\times n$ hyperbolic integral Gram matrix
$$
g=(g_{ij})=(\delta_i, \delta_j), \  \ 1\le i,j \le n,
\tag{1.2.12}
$$ 
or rank $3$, 
$$
\Lambda=(\lambda_1,...\lambda_n)
\tag{1.2.13}
$$
where $\lambda_i \in \bn$, $1\le i \le n$, 
such that the condition \thetag{1.2.4} is valid for all $1\le i \le n-1$ and 
conditions \thetag{1.2.3}, \thetag{1.2.5} --- \thetag{1.2.8} 
are valid for all $1 \le i,j \le n$. The set $(K^n_r)^\prime$ is a subset of 
$K^n_r$ which is characterized by the inequlity 
$$
g_{1n}<-2.
\tag{1.2.14}
$$

To find $K^n_r$, we remark that 
$(g_{ij})$, $1\le i,j \le n-1$, and 
$(\lambda_1,...,\lambda_{n-1})$ should 
belong to $(K^{n-1}_r)^\prime$ and      
$(g_{ij})$, $2\le i,j \le n$, and 
$(\lambda_2,...,\lambda_n)$ as well. Like above, one should calculate 
$g_{1n}=g_{n1}$ using other coefficients, and check that $g_{1n}$ is 
integral non-positive and it satisfies the condition 
\thetag{1.2.7}. Then the pair $(g, \Lambda)$ 
belongs to $K^n_r$. 

The set $K^n_r -(K^n_r)^\prime$ with coprime $\Lambda$ 
gives all $G(A)$ of the size $n$ we are 
looking for.  

The calculations will stop when the set $K^{n+1}_r$ is empty.   

\remark{Remark 1.2.3} For $r=0$ the same algorithm gives  
all possible hyperbolic generalized Cartan matrices of parabolic type 
and of rank $3$ with a lattice Weyl vector 
which are twisted to symmetric generalized Cartan 
matrices. Here one does not need Step 1 and should start with Step 2. 
\endremark

\smallpagebreak 

In Table 2 below we give a realization of the algorithm given above using  
GP/PARI calculator, Version 1.38, by 
C. Batut, D. Bernardi, H. Cohen and M. Olivier. Table 1 gives matrices 
$G(A)$ which we obtain using this program.  
To check (at least) that a matrix $G(A)$ of Table 1 gives a   
generalized Cartan matrix of Theorem 1.2.1, one should construct   
a matrix $C$ of the size $(n+1)\times n$ as follows. 
The first line of $C$  is equal to the first line 
$(\lambda_1,...,\lambda_n)$ 
of the matrix $G(A)$. 
The next lines of $C$ are equal to the lines of   
the symmetric Gram matrix $(\delta_i, \delta_j)$, $1\le i,j\le n$. 
For $i=j$ one has $(\delta_i,\delta_j)=2$. For $i<j$ the 
number 
$-(\delta_i,\delta_j)$ is equal to the coefficient of the $i$-th 
column and $j-i+1$-th line of the matrix  $G(A)$. 
One should check  
that the matrix $C$ has rank $3$; all coefficients $-(\delta_i, \delta_{i+1})$
are not less than $-2$ and 
$\lambda_i(\delta_i,\delta_i)|2\lambda_j(\delta_i,\delta_j)$ for 
all $1\le i,j \le n$.   
 
\vskip20pt  

\centerline{\bf Table 2.} 

\centerline{PARI program for calculation of the 
hyperbolic generallized Cartan} 
\centerline{matrices of elliptic and parablic 
type with a lattice Weyl vector} 
\centerline{which are twisted to symmetric generalized Cartan matrices 
and have rank 3}   

\vskip20pt 

\noindent
$\backslash\backslash$ Unified program for calculation of twisted cases. 
\newline 
$\backslash\backslash$ To calculate with multiplicities up to k, one should  
\newline
$\backslash\backslash$ write main(k). The result will be in the file pari.log  

\smallpagebreak 

\noindent
gram(a,b,c,g)=$\backslash$ 
\newline
g=2*idmat(3);g[1,2]=--a;g[2,1]=--a;g[1,3]=--b;g[3,1]=--b;g[2,3]=--c;g[3,2]=--c;
g;

\smallpagebreak 

\noindent
rho(a,b,c,ta,tb,tc,g)=$\backslash$
\newline 
g=2*idmat(3);g[1,2]=--a;g[2,1]=--a;g[1,3]=--b;g[3,1]=--b;g[2,3]=--c;g[3,2]=--c;
$\backslash$
\newline 
gauss(g,[--ta,--tb,--tc]$^\sim$); 

\smallpagebreak 

\noindent    
rhosq(a,b,c,ta,tb,tc,g,r)=$\backslash$
\newline
g=2*idmat(3);g[1,2]=--a;g[2,1]=--a;g[1,3]=--b;g[3,1]=--b;g[2,3]=--c;g[3,2]=--c;
$\backslash$
\newline
r=gauss(g,[--ta,--tb,--tc]$^\sim$);r$^\sim$*g*r;   

\smallpagebreak 

\noindent
rad(p,rhoo,rho1,a,b,c,ta,tb,tc,g,r,r1,n,k,alpha)=$\backslash$
\newline
n=1;rhoo=vector(n,j,0);rho1=rhoo;$\backslash$
\newline
for(a=0,2,$\backslash$
\newline
for(b=0,14,$\backslash$
\newline
for(c=0,2,$\backslash$
\newline
for(ta=1,p,$\backslash$
\newline
for(tb=1,p,$\backslash$
\newline
for(tc=1,p,$\backslash$
\newline
if(type(tb*a/ta)!=1$||$type(ta*a/tb)!=1$||$type(tb*c/tc)!=1$||$$\backslash$
\newline type(tc*c/tb)!=1$||$type(ta*b/tc)!=1$||$type(tc*b/ta)!=1,,$\backslash$
\newline
g=2*idmat(3);g[1,2]=--a;g[2,1]=--a;g[1,3]=--b;g[3,1]=--b;g[2,3]=--c;g[3,2]=--c;
$\backslash$
\newline 
if(det(g)$>$=0,,$\backslash$
\newline
r1=gauss(g,[--ta,--tb,--tc]$^\sim$);r=r1$^\sim$*g*r1;$\backslash$
\newline if(r$>$=0,,$\backslash$
\newline alpha=1;k=1;$\backslash$ 
\newline while(alpha,$\backslash$
\newline if(rhoo[k]$<$r,k=k+1,$\backslash$
\newline if(rhoo[k]==r,alpha=0,n=n+1;rho1=vector(n,j,0);$\backslash$
\newline for(l=1,k--1,rho1[l]=rhoo[l]);
rho1[k]=r;for(l=k+1,n,rho1[l]=rhoo[l--1]);
$\backslash$
\newline rhoo=rho1;alpha=0))))))))))));rhoo; 

\smallpagebreak 

\noindent
\newline rot(v,s,v1)=$\backslash$
\newline v1=v;for(j=1,s--1,for(k=j+1,s,$\backslash$
\newline if(k$<$s,v1[(j--1)*(2*s--j)/2+k--j]=v[j*(2*s--j--1)/2+k--j],
$\backslash$
\newline v1[(j--1)*(2*s--j)/2+k--j]=v[j])));for(j=1,s,$\backslash$
\newline if(j$<$s,v1[(s--1)*s/2+j]=v[(s--1)*s/2+j+1],
v1[s*(s+1)/2]=v[(s--1)*s/2+1]));v1;

\smallpagebreak 

\noindent
\newline sym(v,s,v1)=,$\backslash$
\newline v1=v;for(j=1,s--1,for(k=j+1,s,$\backslash$
\newline v1[(j--1)*(2*s--j)/2+k--j]=v[(s--k)*(s+k--1)/2+k--j]));$\backslash$
\newline for(j=1,s,$\backslash$
\newline v1[(s--1)*s/2+j]=v[(s--1)*s/2+s--j+1]);v1;

\smallpagebreak 

\noindent 
\newline sorting(a,s,n,m1,m2,a1,a2,a3,t,alpha,beta,v,v1)=$\backslash$
\newline if(type(a)==1$||$type(a)==17,a,$\backslash$
\newline n=matsize(a);m1=n[1];$\backslash$
\newline a1=a;t=1;$\backslash$
\newline v=a1[1,];v1=sym(v,s);if(t$>=$m1,alpha=0,alpha=1);$\backslash$
\newline while(alpha,$\backslash$
\newline a2=matrix(t,n[2],j,k,a1[j,k]);m2=t;$\backslash$
\newline for(j=t+1,m1,$\backslash$
\newline beta=1;$\backslash$
\newline for(k=1,s,v=rot(v,s);v1=rot(v1,s);$\backslash$
\newline if(a1[j,]!=v\&\&a1[j,]!=v1,,beta=0));$\backslash$
\newline if(beta==0,,$\backslash$
\newline m2=m2+1;a3=matrix(m2,n[2],j,k,0);$\backslash$
\newline for(l=1,m2--1,a3[l,]=a2[l,]);a3[m2,]=a1[j,];a2=a3;));$\backslash$
\newline a1=a2;m1=m2;$\backslash$
\newline if(t$>=$m1,alpha=0,$\backslash$
\newline t=t+1;v=vector(n[2],j,a1[t,j]);v1=sym(v,s)));a1);

\smallpagebreak  

\noindent
\newline gr(d,s,b,g)=$\backslash$
\newline g=[2,0,--1;0,2,--2;--1,--2,2];$\backslash$
\newline b=vector((s--1)*s/2,j,0);$\backslash$
\newline for(j=1,s--1,$\backslash$
\newline for(k=j+1,s,$\backslash$
\newline b[(j--1)*(2*s--j)/2+k--j]=--d[j,]*g*d[k,]$^\sim$));b;

\smallpagebreak 

\noindent
$\backslash\backslash$ 
Calculation of the polygonal matrix from the matrix d of vectors 
\newline 
$\backslash\backslash$ delta\_i in bases a,b,c  

\smallpagebreak 

\noindent
polygon(d,s,g,c)=$\backslash$
\newline g=[2,0,--1;0,2,--2;--1,--2,2];$\backslash$
\newline c=matrix(floor(s/2),s,j,k,0);$\backslash$
\newline for(t=1,floor(s/2),for(j=0,s--1,$\backslash$
\newline c[t,j+1]=--d[j+1,]*g*d[lift(mod(j+t,s))+1,]$^\sim$));c;

\smallpagebreak 

\noindent
$\backslash$$\backslash$Calculation of the polygonal matrix from the 
\newline 
$\backslash$$\backslash$gram vector (with multiplicities for delta\_i)

\smallpagebreak 

\noindent
\newline polygon1(b,s,a,c)$\backslash$
\newline a=2*idmat(s);$\backslash$
\newline for(j=1,s--1,for(k=j+1,s,$\backslash$
\newline a[j,k]=--b[(j--1)*(2*s--j)/2+k--j];a[k,j]=a[j,k]));$\backslash$
\newline c=matrix(floor(s/2)+2,s,j,k,0);$\backslash$
\newline for(t=1,floor(s/2),for(j=0,s--1,$\backslash$
\newline c[t+2,j+1]=--a[j+1,lift(mod(j+t,s))+1]));$\backslash$
\newline for(j=1,s,c[1,j]=b[s*(s--1)/2+j]);c; 

\smallpagebreak 

\noindent
$\backslash$$\backslash$calculation of the symmetric 
generalized Cartan 
\newline $\backslash\backslash$ 
matrix (Gram matrix of m\_i$\backslash$delta\_i) from the Gram vector

\smallpagebreak 

\noindent 
symcartan(b,s,a,sa)=$\backslash$
\newline a=2*idmat(s);$\backslash$
\newline for(j=1,s--1,for(k=j+1,s,$\backslash$
\newline a[j,k]=--b[(j--1)*(2*s--j)/2+k--j];a[k,j]=a[j,k]));$\backslash$
\newline sa=matrix(s,s,j,k,b[s*(s--1)/2+j]*b[s*(s--1)/2+k]*a[j,k]);sa;

\smallpagebreak 

\noindent
$\backslash\backslash$ function eq3(p,r)
\newline eq3(p,r,s,m,u,r1,alpha)=s=matrix(60,6,j,k,0);m=0;s[1,]=[0,0,0,0,0,0];
$\backslash$
\newline for(a=0,2,for(c=0,2,for(ta=1,p,for(tc=1,p,for(tb=1,p,$\backslash$
\newline b=1;alpha=1;$\backslash$
\newline while(alpha,$\backslash$
\newline if(det(gram(a,b,c))$>$=0,b=b+1,$\backslash$
\newline r1=rhosq(a,b,c,ta,tb,tc);$\backslash$
\newline if(r1$<$r,b=b+1,);$\backslash$
\newline if(r1$>$r,alpha=0,);$\backslash$
\newline if(r1==r,$\backslash$
\newline if(type(tb*a/ta)==1,$\backslash$
\newline if(type(ta*a/tb)==1,$\backslash$
\newline if(type(tb*c/tc)==1,$\backslash$
\newline if(type(tc*c/tb)==1,$\backslash$
\newline if(type(ta*b/tc)==1,$\backslash$
\newline if(type(tc*b/ta)==1,$\backslash$
\newline m=m+1;s[m,]=[a,b,c,ta,tb,tc],),),),),),);$\backslash$
\newline alpha=0,))$\backslash$
\newline ))))));$\backslash$
\newline if(m==0,0,u=matrix(m,6,j,k,0);for(j=1,m,u[j,]=s[j,]);u);  
 
\smallpagebreak 

\noindent 
\newline $\backslash\backslash$ function eq4(p,r),s=4
\newline eq4(p,r,a,aa,aaa,n,m,u,v,w,x,y,z,
talpha,tbeta,tgamma,tdelta,ra)=$\backslash$
\newline u=eq3(p,r);$\backslash$
\newline if(type(u)==1,0,$\backslash$
\newline n=matsize(u);n=n[1];v1=matrix(n,6,j,k,0);m=0;$\backslash$
\newline for(x=1,n,if(u[x,2]$>$2,m=m+1;v1[m,]=u[x,],));$\backslash$
\newline if(m==0,a=0,$\backslash$
\newline a=matrix(m,6,j,k,0);for(x=1,m,a[x,]=v1[x,])));$\backslash$
\newline if(type(a)==1,0,$\backslash$
\newline n=matsize(a);n=n[1];aa=matrix(100,10,j,k,0);m=0;$\backslash$
\newline for(j=1,n,u=a[j,1];v=a[j,2];w=a[j,3];talpha=a[j,4];$\backslash$
\newline tbeta=a[j,5];tgamma=a[j,6];$\backslash$
\newline for(k=1,n,$\backslash$
\newline if([w,tbeta,tgamma]!=[a[k,1],a[k,4],a[k,5]],,$\backslash$
\newline y=a[k,2];z=a[k,3];
tdelta=a[k,6];ra=rho(u,v,w,talpha,tbeta,tgamma);$\backslash$
\newline x=(tdelta--ra[2]*y--ra[3]*z)/ra[1];$\backslash$
\newline if(x$<$0,,$\backslash$
\newline if(type(x)!=1$||$type(tdelta*x/talpha)!=1$||$
type(talpha*x/tdelta)!=1,,$\backslash$
\newline if(det([2,--u,--v,--x;--u,2,--w,--y;--v,--w,2,--z;--x,--y,--z,2])!=0,,
$\backslash$
\newline m=m+1;aa[m,]=[u,v,x,w,y,z,talpha,tbeta,tgamma,tdelta]))))));
$\backslash$
\newline if(m==0,aaa=0,aaa=matrix(m,10,j,k,0);
for(j=1,m,aaa[j,]=aa[j,]));$\backslash$
aaa);

\smallpagebreak 

\noindent
eq(p,r,a,aa,n,s,alpha,m,mm,vv,ww,u,v,w,u1,v1,w1,g,g1,$\backslash$
\newline delta4,delta4du,deltalast,deltalastdu,deltalast1)=$\backslash$
\newline s=3;a=eq3(p,r);$\backslash$
\newline if(type(a)==1,,n=matsize(a);$\backslash$
\newline aa1=matrix(n[1],n[2],j,k,0);m1=0;$\backslash$
\newline for(j=1,n[1],$\backslash$
\newline if(a[j,s--1]$>$2,,$\backslash$
\newline aaa=vector(s,k,a[j,n[2]--s+k]);$\backslash$
\newline if(content(aaa)!=1,,m1=m1+1;aa1[m1,]=a[j,])));$\backslash$
\newline if(m1==0,a1=0,$\backslash$
\newline a1=matrix(m1,n[2],j,k,aa1[j,k]));$\backslash$
\newline a1=sorting(a1,s);$\backslash$
\newline if(type(a1)==1,,n1=matsize(a1);for(j=1,n1[1],$\backslash$
\newline qq=a1[j,];$\backslash$
\newline pprint(qq);texprint(qq);$\backslash$
\newline pprint(symcartan(qq,s));texprint(symcartan(qq,s));$\backslash$
\newline pprint(polygon1(qq,s));texprint(polygon1(qq,s)))));$\backslash$
\newline s=4;a=eq4(p,r);$\backslash$
\newline alpha=1;$\backslash$
\newline while(alpha,$\backslash$
\newline if(type(a)==1,alpha=0,n=matsize(a);$\backslash$
\newline aa=matrix(n[1],n[2],j,k,0);m=0;$\backslash$
\newline aa1=matrix(n[1],n[2],j,k,0);m1=0;$\backslash$
\newline for(j=1,n[1],$\backslash$
\newline if(a[j,s--1]$>$2,m=m+1;aa[m,]=a[j,],$\backslash$
\newline aaa=vector(s,k,a[j,n[2]--s+k]);$\backslash$
\newline if(content(aaa)!=1,,m1=m1+1;aa1[m1,]=a[j,])));$\backslash$
\newline if(m1==0,a1=0,$\backslash$
\newline a1=matrix(m1,n[2],j,k,aa1[j,k]));$\backslash$
\newline a1=sorting(a1,s);$\backslash$
\newline if(type(a1)==1,,n1=matsize(a1);for(j=1,n1[1],$\backslash$
\newline qq=a1[j,];$\backslash$
\newline pprint(qq);texprint(qq);$\backslash$
\newline pprint(symcartan(qq,s));texprint(symcartan(qq,s));$\backslash$
\newline pprint(polygon1(qq,s));texprint(polygon1(qq,s))));$\backslash$
\newline if(m==0,a=0;alpha=0,$\backslash$
\newline a=matrix(m,n[2],j,k,aa[j,k]);$\backslash$
\newline aa=matrix(100,n[2]+s+1,j,k,0);mm=0;$\backslash$
\newline for(j=1,m,for(k=1,m,$\backslash$
\newline vv=vector(n[2]--s,l,0);$\backslash$
\newline for(l=1,n[2]--2*s+1,vv[l]=a[j,s--1+l]);$\backslash$
\newline for(l=n[2]--2*s+2,n[2]--s,vv[l]=a[j,s+l]);$\backslash$
\newline ww=vector(n[2]--s,l,0);$\backslash$
\newline for(y=0,s--3,$\backslash$
\newline for(l=y*(2*s--y--3)/2+1,
y*(2*s--y--3)/2+s--2--y,ww[l]=a[k,y+l]));$\backslash$
\newline for(l=n[2]--2*s+2,n[2]--s,ww[l]=a[k,l+s--1]);$\backslash$
\newline if(vv!=ww,,$\backslash$
\newline u=a[j,1];v=a[j,2];w=a[j,s];g=gram(u,v,w);$\backslash$
\newline delta4du=[--a[j,3],--a[j,s+1],--a[j,2*s--2]];$\backslash$
\newline delta4=delta4du*g$^\wedge$--1;$\backslash$
\newline c=[0,1,0;0,0,1;delta4[1],delta4[2],delta4[3]];$\backslash$
\newline u1=a[k,1];v1=a[k,2];w1=a[k,s];g1=gram(u1,v1,w1);$\backslash$
\newline deltalastdu=[--a[k,s--1],--a[k,s--1+s--2],--a[k,s--1+s--2+s--3]];
$\backslash$
\newline deltalast1=deltalastdu*g1$^\wedge$--1;$\backslash$
\newline deltalast=deltalast1*c;$\backslash$
\newline x=--deltalast*g*[1,0,0]$^\sim$;$\backslash$
\newline if(x$<$0$||$type(x)!=1$||$type(x*a[j,n[2]--s+1]/a[k,n[2]])!=1$||$
$\backslash$
\newline type(x*a[k,n[2]]/a[j,n[2]--s+1])!=1,,$\backslash$
\newline mm=mm+1;$\backslash$
\newline for(l=1,s--1,aa[mm,l]=a[j,l]);aa[mm,s]=x;$\backslash$
\newline for(l=s+1,s*(s+1)/2,aa[mm,l]=a[k,l--s]);$\backslash$
\newline aa[mm,s*(s+1)/2+1]=a[j,(s--1)*s/2+1];$\backslash$
\newline for(l=s*(s+1)/2+2,(s+1)*(s+2)/2,aa[mm,l]=a[k,l--s--1])))));
$\backslash$
\newline if(mm==0,alpha=0;a=0,$\backslash$
\newline a=matrix(mm,(s+1)*(s+2)/2,l1,l2,aa[l1,l2]);s=s+1)))); 

\smallpagebreak 

\noindent
$\backslash$l;

\smallpagebreak 

\noindent 
main(p,n)=$\backslash$
\newline radius=rad(p);n=matsize(radius);$\backslash$
\newline n=n[2];print("n=",n);$\backslash$
\newline for(j=1,n--1,r=radius[j];print("j=",j,"; r=",r);eq(p,r));

\enddemo

\remark{Remark 1.2.4} Table 1 contains $60$ matrices $G(A)$. 
Seven of them correspond to compact case (i.e., a fundamental polygon $\M$ 
for $W$ is compact in the hyperbolic plane, it has only finite vertices).  
These cases correspond to matrices $G(A)$ with second line 
without $2$ (for $2$ the corresponding angle 
of $\M$ is equal to $0$, for $1$ 
it is equal to $\pi/3$, for $0$ it is equal to $\pi/2$). 

Between $7$ compact cases there are $4$ non-twisted ones (i.e. they 
give symmetric generalized Cartan matrices). Non-twisted cases correspond 
to matrices $G(A)$ with first line containing $1$ only. 
Between $53$ non-compact cases there are $12$ non-twisted ones. 

We consider these last $12$ cases below.
\endremark

\subhead
1.3. Symmetric hyperbolic generalized Cartan matrices of rank $3$ 
having elliptic type and a lattice Weyl vector. Non-compact case 
\endsubhead  

From Theorem 1.2.1, we get $12$ 
symmetric generalized Cartan matrices $A$ corresponding to $12$ non-compact 
and non-twisted matrices $G(A)$ of Table 1 (see Remark 1.2.4). We announced 
this list in \cite{GN4}.   

\proclaim{Theorem 1.3.1} There are exactly 12 
symmetric hyperbolic generalized Cartan matrices of elliptic type and 
of rank $3$ which have Weyl group with a non-compact (i.e.
with an infinite vertex)
fundamental polygon and have a lattice Weyl vector.
They are matrices $A_{1,0}$ -- $A_{3,III}$ below:
$$
\align
A_{1,0}&=
\pmatrix
\hphantom{-}{2}&\hphantom{-}{0}&{-1}\cr
\hphantom{-}{0}&\hphantom{-}{2}&{-2}\cr
{-1}&{-2}&\hphantom{-}{2}\cr
\endpmatrix,\qquad
A_{1,I}=
\pmatrix
\hphantom{-}{2}&{-2}&{-1}\cr
{-2}&\hphantom{-}{2}&{-1}\cr
{-1}&{-1}&\hphantom{-}{2}\cr
\endpmatrix,\\
A_{1,II}&=
\pmatrix
\hphantom{-}{2}&{-2}&{-2}\cr
{-2}&\hphantom{-}{2}&{-2}\cr
{-2}&{-2}&\hphantom{-}{2}\cr
\endpmatrix,\quad\,
A_{1,III}=
\pmatrix
\hphantom{-}{2}&{-2}&{-6}&{-6}&{-2}\cr
{-2}&\hphantom{-}{2}&\hphantom{-}{0}&{-6}&{-7}\cr
{-6}&\hphantom{-}{0}&\hphantom{-}{2}&{-2}&{-6}\cr
{-6}&{-6}&{-2}&\hphantom{-}{2}&{0}\cr
{-2}&{-7}&{-6}&\hphantom{-}{0}&\hphantom{-}{2}\cr
\endpmatrix ,
\endalign
$$
$$
A_{2,0}=
\pmatrix
\hphantom{-}{2}&{-2}&{-2}\cr
{-2}&\hphantom{-}{2}&\hphantom{-}{0}\cr
{-2}&\hphantom{-}{0}&\hphantom{-}{2}\cr
\endpmatrix,
$$
$$
A_{2,I}=
\pmatrix
\hphantom{-}{2}&{-2}&{-4}&\hphantom{-}{0}\cr
{-2}&\hphantom{-}{2}&\hphantom{-}{0}&{-4}\cr
{-4}&\hphantom{-}{0}&\hphantom{-}{2}&{-2}\cr
\hphantom{-}{0}&{-4}&{-2}&\hphantom{-}{2}\cr
\endpmatrix ,
\qquad
A_{2,II}=
\pmatrix
\hphantom{-}{2}&{-2}&{-6}&{-2}\cr
{-2}&\hphantom{-}{2}&{-2}&{-6}\cr
{-6}&{-2}&\hphantom{-}{2}&{-2}\cr
{-2}&{-6}&{-2}&\hphantom{-}{2}\cr
\endpmatrix ,
$$
$$
A_{2,III}=
\pmatrix
\hphantom{-}{2}&{-2}&{-8}&{-16}&{-18}&{-14}&{-8}&\hphantom{-}{0}\cr
{-2}&\hphantom{-}{2}&\hphantom{-}{0}&{-8}&{-14}&{-18}&{-16}&{-8}\cr
{-8}&\hphantom{-}{0}&\hphantom{-}{2}&{-2}&{-8}&{-16}&{-18}&{-14}\cr
{-16}&{-8}&{-2}&\hphantom{-}{2}&\hphantom{-}{0}&{-8}&{-14}&{-18}\cr
{-18}&{-14}&{-8}&\hphantom{-}{0}&\hphantom{-}{2}&{-2}&{-8}&{-16}\cr
{-14}&{-18}&{-16}&{-8}&{-2}&\hphantom{-}{2}&\hphantom{-}{0}&{-8}\cr
{-8}&{-16}&{-18}&{-14}&{-8}&\hphantom{-}{0}&\hphantom{-}{2}&{-2}\cr
\hphantom{-}{0}&{-8}&{-14}&{-18}&{-16}&{-8}&{-2}&\hphantom{-}{2}\cr
\endpmatrix ,
$$
$$
A_{3,0}=
\pmatrix
\hphantom{-}{2}&{-2}&{-2}\cr
{-2}&\hphantom{-}{2}&{-1}\cr
{-2}&{-1}&\hphantom{-}{2}\cr
\endpmatrix ,
\qquad\qquad\ A_{3,I}=
\pmatrix
\hphantom{-}{2}&{-2}&{-5}&{-1}\cr
{-2}&\hphantom{-}{2}&{-1}&{-5}\cr
{-5}&{-1}&\hphantom{-}{2}&{-2}\cr
{-1}&{-5}&{-2}&\hphantom{-}{2}\cr
\endpmatrix ,
$$
$$
A_{3,II}=
\pmatrix
\hphantom{-}{2}&{-2}&{-10}&{-14}&{-10}&{-2}\cr
{-2}&\hphantom{-}{2}&{-2}&{-10}&{-14}&{-10}\cr
{-10}&{-2}&\hphantom{-}{2}&{-2}&{-10}&{-14}\cr
{-14}&{-10}&{-2}&\hphantom{-}{2}&{-2}&{-10}\cr
{-10}&{-14}&{-10}&{-2}&\hphantom{-}{2}&{-2}\cr
{-2}&{-10}&{-14}&{-10}&{-2}&\hphantom{-}{2}\cr
\endpmatrix ,
$$
$$
A_{3,III}=
$$
$$
\pmatrix
\hphantom{-}{2}&{-2}&{-11}&{-25}&{-37}&{-47}&{-50}
&{-46}&{-37}&{-23}&{-11}&{-1}\cr
{-2}&\hphantom{-}{2}&{-1}&{-11}&{-23}&{-37}&
{-46}&{-50}&{-47}&{-37}&{-25}&{-11}\cr
{-11}&{-1}&\hphantom{-}{2}&{-2}&{-11}&{-25}
&{-37}&{-47}&{-50}&{-46}&{-37}&{-23}\cr
{-25}&{-11}&{-2}&\hphantom{-}{2}&{-1}&{-11}
&{-23}&{-37}&{-46}&{-50}&{-47}&{-37}\cr
{-37}&{-23}&{-11}&{-1}&\hphantom{-}{2}&{-2}
&{-11}&{-25}&{-37}&{-47}&{-50}&{-46}\cr
{-47}&{-37}&{-25}&{-11}&{-2}&\hphantom{-}{2}
&{-1}&{-11}&{-23}&{-37}&{-46}&{-50}\cr
{-50}&{-46}&{-37}&{-23}&{-11}&{-1}&\hphantom{-}{2}
&{-2}&{-11}&{-25}&{-37}&{-47}\cr
{-46}&{-50}&{-47}&{-37}&{-25}&{-11}&{-2}
&\hphantom{-}{2}&{-1}&{-11}&{-23}&{-37}\cr
{-37}&{-47}&{-50}&{-46}&{-37}&{-23}&{-11}&{-1}
&\hphantom{-}{2}&{-2}&{-11}&{-25}\cr
{-23}&{-37}&{-46}&{-50}&{-47}&{-37}
&{-25}&{-11}&{-2}&\hphantom{-}{2}&{-1}&{-11}\cr
{-11}&{-25}&{-37}&{-47}&{-50}&{-46}&{-37}
&{-23}&{-11}&{-1}&\hphantom{-}{2}&{-2}\cr
{-1}&{-11}&{-23}&{-37}&{-46}&{-50}
&{-47}&{-37}&{-25}&{-11}&{-2}&\hphantom{-}{2}\cr
\endpmatrix .
$$
\endproclaim

\demo{Proof} We leave the trivial calculation of matrices $A$ 
from matrices $G(A)$ to a reader. In Table 1, non-compact non-twisted 
cases give respectively: 
\newline  
$r=-23/2$ gives $A_{1,0}$; 
\newline
$r=-4$ gives $A_{1,I}$; 
\newline
$r=-7/2$ gives $A_{2,0}$;
\newline
$r=-13/6$ gives $A_{3,0}$;
\newline
$r=-3/2$ gives $A_{1,II}$;
\newline
$r=-1$ gives $A_{2,I}$;
\newline
$r=-2/3$ gives $A_{3,I}$;
\newline
$r=-1/2$ gives $A_{2,II}$;
\newline
$r=-7/18$ gives $A_{1,III}$;
\newline
$r=-1/6$ gives $A_{3,II}$; 
\newline
$r=-1/8$ gives $A_{2,III}$;
\newline
$r=-1/24$ gives $A_{3,III}$. 
\enddemo

We numerate $12$ matrices of Theorem 1.3.1 by 
two indices $i=1,2,3$ and $j=0,I,II, III$. 
For non-twisted cases the matrices $A_{i,j}$ are the 
Gram matrices of 
sets $P(\M_{i,j})$ of the orthogonal vectors (normalized by square $2$)  
to fundamental polygons $\M_{i,j}$ of reflection groups $W_{i,j}$ 
on a hyperbolic plane.  
The polygons $\M_{i,\ast }$ corresponding to $A_{i,\ast }$, 
$i=1,2,3$,   
are naturally composed from a triangle with two angles $0$ and $\pi/2$ and 
the third angle $\pi/3$, $\pi/4$ and $\pi/6$ respectively. The 
second index $j$ 
shows how we get the polygon $\M_{\ast ,j}$ from this triangle.  
Below we describe the polygons $\M_{i,j}$, the sets $P(\M_{i,j})$,  
the lattices $M_{i,j}$ and the lattice Weyl vectors $\rho_{i,j}$ 
corresponding to $A_{i,j}$.    

\vskip20pt 
We use the following notations for lattices and their 
discriminant forms: 
\newline
$\langle A \rangle$ is a lattice with a matrix $A$; 
\newline
$\oplus$ denotes an orthogonal sum of lattices, 
$K^n$ is orthogonal sum of $n$ copies of a lattice 
$K$ (the same for forms);
\newline 
$K(n)$ denotes a lattice which 
is obtained multiplying by $n\in \bq$ the form of a lattice 
$K$;
\newline
$nK$, $n \in \bq$, denotes a lattice $(nM,\phi|_{nM})$ if $K=(M,\phi)$, 
evidently, $nK\cong K(n^2)$;   
\newline
$b_K$ and $q_K$ are the discriminant bilinear and  
quadratic forms respectively on the discriminant group $A_K=K^\ast/K$ 
of a lattice $K$; for discriminant forms 
we use notations from \cite{N2}; 
\newline
 $b_\epsilon^{(p)}(p^k)$, $q_\epsilon^{(p)}(p^k)$ are the  
discriminant bilinear and quadratic form respectively of a one-dimensional 
$p$-adic lattice (i.e. over the ring $\bz_p$ of $p$-adic 
integers)  
$K_\epsilon^{(p)}(p^k)=\langle \epsilon p^k\rangle$, 
$\epsilon \in \bz_p^\ast$; 
\newline
$u_-^{(2)}(2^k)$, $v_-^{(2)}(2^k)$ and $u_+^{(2)}(2^k)$,  
$v_+^{(2)}(2^k)$ are discriminant bilinear and quadratic forms of $2$-adic  
lattices $U^{(2)}(2^k)=U(2^k)\otimes \bz_2$ and 
$V^{(2)}(2^k)=A_2\otimes \bz_2$ respectively 
where 
$$
U=\langle
\matrix
0&-1\\
-1&0
\endmatrix 
\rangle, 
\ \ \ \ 
A_2=
\langle
\matrix
2&-1\\
-1&2
\endmatrix 
\rangle; 
$$
\newline
$W^{(m_1,...,m_k)}(K)$ denotes a subgroup of $O^+(K)$ 
of a lattice $K$ generated by reflections $s_\delta$   
in all primitive elements $\delta \in K$ with square 
$(\delta,\delta) \in \{m_1,...,m_k\}$; here 
$O^+(K)$ is a subgroup of the full orthogonal group $O(K)$ which 
fixes the cone $V^+(K)$.   
   
\subsubhead
1.3.1. Generalized Cartan matrices $A_{1,j}$
\endsubsubhead

In Figure 1 we consider the following construction on a hyperbolic 
plane. Let $ABC$ be a triangle with the angles 
$0$, $\pi/3$ and $\pi/2$ respectively. Considering the image of 
$ABC$ by the reflection in $AC$, we get a triangle $BB_1A$ with 
the angles $\pi/3$, $\pi/3$ and $0$ respectively. Considering the image of 
the triangle $ABC$ by the group (of type $D_3$) generated by reflections 
in $BC$ and $BA$, we get a right triangle $AA_2A_3$ (here $AA_2$ contains 
$C$) with zero angles. The right triangle $AA_2A_3$ is divided 
by medians $AC_2$, $A_2C_3$ and $A_3C$ (all containing the center $B$).
Considering a reflection in $A_2A_3$ of the triangle $A_2A_3C_3$ and 
a reflection in $AA_3$ of the triangle $AC_2A_3$, we get a pentagon 
$AA_2C_3^\prime A_3C_2^\prime$ with the 
angles $0$, $0$, $\pi/2$, $0$, $\pi/2$ 
respectively. Here $C_3^\prime$ is the image of $C_3$ and $C_2^\prime$ 
is the image of $C_2$. 

It is easy to see that the matrix $A_{1,0}$ is the Gram matrix (of 
orthogonal vectors with square $2$ to the 
sides) of the triangle $\M_{1,0}=ABC$. 
The matrix $A_{1,I}$ is the Gram matrix of the triangle $\M_{1,I}=BB_1A$, 
the matrix $A_{1,II}$ is the Gram matrix of the triangle 
$\M_{1,II}=AA_2A_3$, and 
$A_{1,III}$ is the Gram matrix 
of the pentagon $\M_{1,III}=AA_2C_3^\prime A_3C_2^\prime$.

\midinsert
$$\vbox{\centerline{\epsfxsize=5in\epsfbox{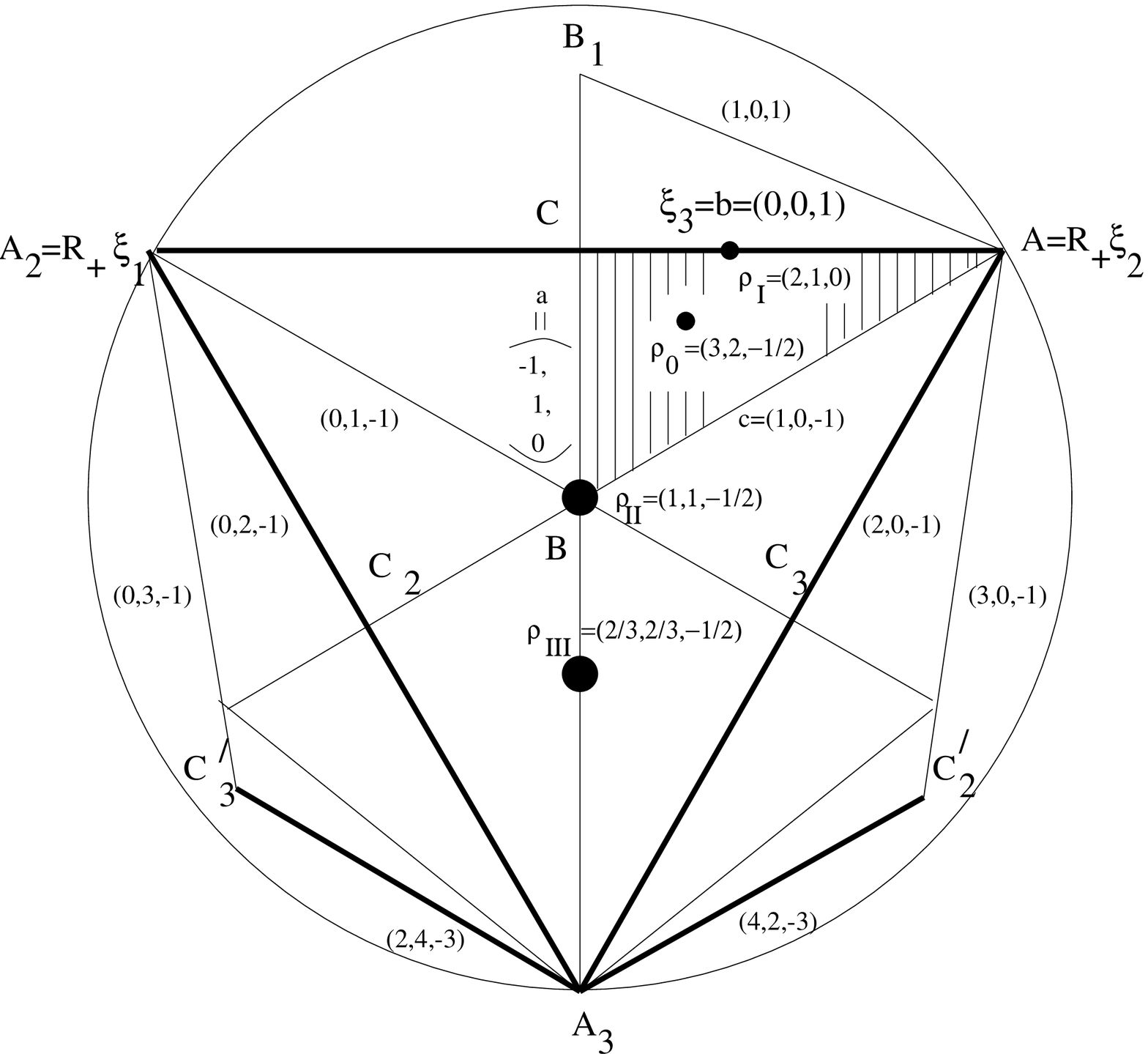}}
\vskip20pt 
\centerline{Figure 1. Fundamental polygons of $A_{1,j}$}}$$
\endinsert

We introduce a hyperbolic lattice 
$M_{1,0}\cong U\oplus \langle 2 \rangle $ with the corresponding 
standard basis    
$\xi_1$, $\xi_2$, $\xi_3$ having the Gram matrix 
$$
\left(
\matrix
0 &-1&0\\
-1&0&0\\
0&0& 2
\endmatrix\right). 
$$
Thus, $\xi_1, \xi_2$ and $\xi_3$ are the standard basis of orthogonal 
lattices $U$ and $\langle 2 \rangle$ respectively.
We identify the hyperbolic planes of the lattice $M_{1,0}$ and 
Figure 1 as follows. 
Point $A_2=\br_{++}\xi_1$, $A=\br_{++}\xi_2$, elements  
$a=-\xi_1+\xi_2$, $b=\xi_3$ and $c=\xi_1-\xi_3$ are orthogonal 
to sides $CB$, $AC$ and $AB$ of the triangle $ABC$ and 
are directed outward. 

The triangle $\M_{1,0}=BCA$ has    
$$
P(\M_{1,0})=\{\delta_1=a=-\xi_1+\xi_2,\  
\delta_2=b=\xi_3, \ \delta_3=c=\xi_1-\xi_3\}  
\tag{1.3.1}
$$
with the Gram matrix $A_{1,0}$ and the lattice Weyl vector 
$$
\rho_{1,0}=3\xi_1+2\xi_2-(1/2)\xi_3.
\tag{1.3.2}
$$    
The lattice $M_{1,0}$ is generated by $P(\M_{1,0})$ and is then the 
lattice of $A_{1,0}$. Thus, 
$$
M_{1,0}\cong U\oplus \langle 2 \rangle.  
\tag{1.3.3}
$$  
We have 
$$
W_{1,0}=W^{(2)}(M_{1,0}), \ \ \ \Sym(P(\M_{1,0}))\  \text{is trivial}.  
\tag{1.3.4}
$$
The \thetag{1.3.4} are  well-known. For example, they follow from 
Vinberg's algorithm \cite{V1} for calculation of a fundamental polyhedron 
of a hyperbolic refection group.  

The triangle $\M_{1,I}=B_1AB$ has 
$$
P(\M_{1,I})=\{\delta_1=\xi_1+\xi_3,\  
\delta_2=\xi_1-\xi_3,\  \delta_3=-\xi_1+\xi_2\} 
\tag{1.3.5}
$$
with the Gram matrix $A_{1,I}$ and the lattice Weyl vector 
$$
\rho_{1,I}=2\xi_1+\xi_2.
\tag{1.3.6}
$$ 
The lattice $M_{1,I}=[a,2b,c]$ 
is a sublattice of the lattice $M_{1,I}$. 
It has the discriminant form 
$$
q_{M_{1,I}}=q^{(2)}_5(8).
\tag{1.3.7} 
$$  
We have  
$$
W_{1,I}=W^{(2)}(M_{1,I}),\ \ \Sym(\M_{1,I})=D_1 
\tag{1.3.8}
$$ 
where the second group 
is generated by the reflection in $2\xi_3=\delta_1-\delta_2$. Here 
and in similar calculations below we use results of \cite{N7} where 
all 2-reflective integral hyperbolic lattices $M$ of rank $3$ were 
classified and their reflection groups $W^{(2)}(M)$ were calculated. 
A lattice $M$ is called {\it $2$-reflective}  
if $W^{(2)}(M)$ has finite index in $O(M)$.  

The triangle $\M_{1,II}=AA_3A_2$ has 
$$
P(\M_{1,II})=\{\delta_1=2\xi_1-\xi_3,\  \delta_2=2\xi_2-\xi_3,\  
\delta_3=\xi_3\}
\tag{1.3.9}
$$
with the Gram matrix $A_{1,II}$ and the lattice Weyl vector 
$$
\rho_{1,II}=\xi_1+\xi_2-(1/2)\xi_3.
\tag{1.3.10}
$$
 The lattice 
$M_{1,II}=[2\xi_1, 2\xi_2, \xi_3]=[2a,b,2c]$ is 
$$
M_{1,II}\cong U(4)\oplus \langle 2 \rangle.
\tag{1.3.11}
$$ 
We have  
$$
W_{1,II}=W^{(2)}(M_{1,II}),\ \ \ \Sym(P(\M_{1,II}))=D_3
\tag{1.3.12}
$$
where the second group 
is generated by the reflections in $c=\xi_1-\xi_3$ and 
$a=-\xi_1+\xi_2$.

The pentagon $\M_{1,III}=AA_2C_3^\prime A_3C_2^\prime$ has 
$$
P(\M_{1,III})=\{\delta_1=\xi_3,\  \delta_2=3\xi_2-\xi_3,  
\delta_3=2\xi_1+4\xi_2-3\xi_3, \delta_4=4\xi_1+2\xi_2-3\xi_3, 
\delta_5=3\xi_1-\xi_3\}
\tag{1.3.13}
$$
with the Gram matrix $A_{1,III}$ and the lattice Weyl vector 
$$
\rho_{1,III}=(2/3)\xi_1+(2/3)\xi_2-(1/2)\xi_3.
\tag{1.3.14}
$$
The lattice $M_{1,III}=[a,b,3c]$ 
has the discriminant form 
$$
q_{M_{1,III}}=q_1^{(2)}(2)\oplus q_1^{(3)}(9).
\tag{1.3.15}
$$
We have 
$$
W_{1,III}\triangleleft W^{(2)}(M_{1,III})\ \ \text{and}\ \ 
W^{(2)}(M_{1,III})/W_{1,III}=D_1
\tag{1.3.16}
$$
where the second group is generated by the reflection $s_a$. 
The reflection $s_a$ also generates 
the group $\Sym(P(\M_{1,III}))$. 
 
We remark that all lattices $M_{1,0}$, $M_{1,I}$, $M_{1,II}$, 
$M_{1,III}$ are characterized by their discriminant quadratic forms:  
i.e. a lattice of signature $(2,1)$ with one of 
these discriminant quadratic forms is unique.  

\subsubhead
1.3.2. Generalized Cartan matrices $A_{2,j}$
\endsubsubhead

\midinsert
$$\vbox{\centerline{\epsfxsize=5in\epsfbox{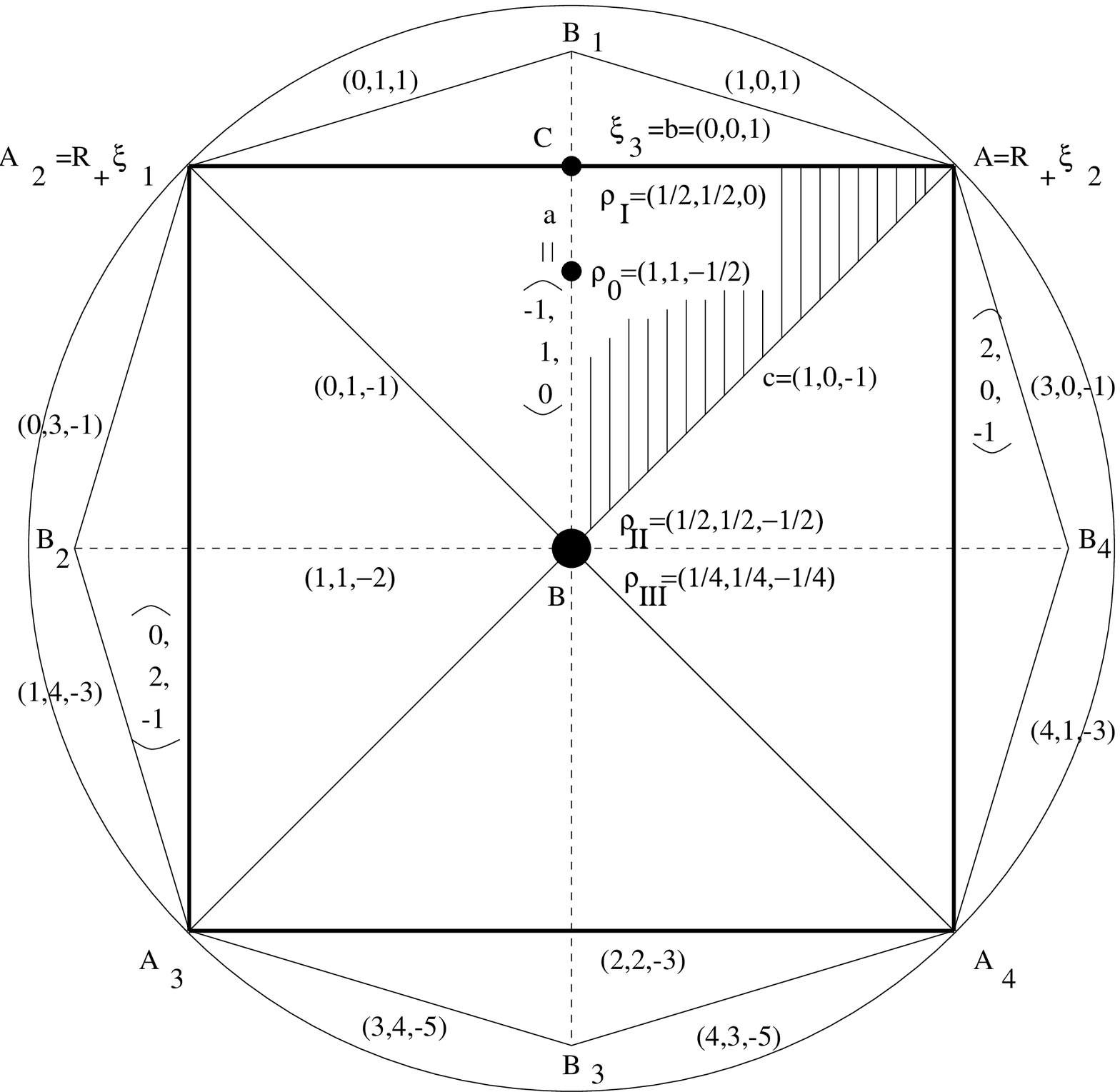}}
\vskip20pt 
\centerline{Figure 2. Fundamental polygons of $A_{2,j}$}}$$
\endinsert

On Figure 2 $ABC$ is a triangle with the angles 
$0$, $\pi/4$ and $\pi/2$ respectively. Considering the orbit of 
$ABC$ by the group generated by the reflection in $BC$, 
we get a triangle $ABA_2$ with the angles $0$, $\pi/2$ and 
$0$ respectively. Considering the orbit of the triangle 
$ABC$ by the group (of type $D_2$) generated by the 
reflections in $BC$ and 
$AC$, we get a quadrangle $ABA_2B_1$ with 
the angles $0$, $\pi/2$, $0$ and $\pi/2$ respectively. 
Considering the orbit of 
the triangle $ABC$ by the group (of type $D_4$) generated by 
the reflections 
in $BC$ and $BA$, we get a right quadrangle $AA_2A_3A_4$ 
(here $AA_2$ contains $C$) with zero angles and the center $B$. 
Considering images of the center $B$ by reflections in the sides 
$AA_2$, $A_2A_3$, $A_3A_4$ and $A_4A$, we get points $B_1$, $B_2$, $B_3$ 
and $B_4$ respectively which together with $A$, $A_2$, $A_3$ and 
$A_4$ define an $8$-gon $AB_1A_2B_2A_3B_3A_4B_4$ with the 
angles $0$, $\pi/2$, 
$0$, $\pi/2$, $0$, $\pi/2$,  $0$, $\pi/2$ respectively. This 
$8$-gon also is the orbit of the triangle $ABB_1$ with the angles 
$0$, $\pi/4$, $\pi/4$ 
respectively by the group (of type $D_4$) generated by reflections 
in his sides $BA$, $BB_1$. 

It is easy to see that the matrix $A_{2,0}$ is the Gram matrix (of 
orthogonal vectors with square $2$ to the sides) of the triangle 
$\M_{2,0}=ABA_2$. 
The matrix $A_{2,I}$ is the Gram matrix of the quadrangle 
$\M_{2,I}=ABA_2B_1$.   
The matrix $A_{2,II}$ is the Gram matrix of the right quadrangle  
$\M_{2,II}=AA_2A_3A_4$, and 
$A_{2,III}$ is the Gram matrix 
of the $8$-gon $\M_{2,III}=AB_1A_2B_2A_3B_3A_4B_4$.    
 
 We introduce a hyperbolic lattice 
$M_{2,0}\cong U(2)\oplus \langle 2 \rangle $ with the corresponding 
standard bases    
$\xi_1$, $\xi_2$, $\xi_3$ having the Gram matrix 
$$
\left(
\matrix
0 &-2&0\\
-2&0&0\\
0&0& 2
\endmatrix\right). 
$$
We identify the hyperbolic plane of the lattice $M_{2,0}$ and 
of Figure 2:  
Point $A_2=\br_{++}\xi_1$, $A=\br_{++}\xi_2$, elements  
$a=-\xi_1+\xi_2$, $b=\xi_3$ and $c=\xi_1-\xi_3$ are orthogonal 
to sides $CB$, $AC$ and $AB$ of the triangle $ABC$ and 
directed outward. Here $(a,a)=4$ and $(b,b)=(c,c)=2$. 

We have for $\M_{2,0}=AA_2B$   
$$
P(\M_{2,0})=\{\delta_1=\xi_3,\  
\delta_2=\xi_2-\xi_3,\ \delta_3=\xi_1-\xi_3\}  
\tag{1.3.17}
$$
with the Gram matrix $A_{2,0}$ and the lattice Weyl vector 
$$
\rho_{2,0}=\xi_1+\xi_2-(1/2)\xi_3.
\tag{1.3.18}
$$    
The lattice $M_{2,0}$ generated by $P(\M_{2,0})$ is then $M_{2,0}=[a,b,c]$, 
and   
$$
M_{2,0}\cong U(2)\oplus \langle 2 \rangle.  
\tag{1.3.19}
$$  
We have 
$$
W_{2,0}=W^{(2)}(M_{2,0}), \ \ \ \Sym(P(\M_{2,0}))=D_1    
\tag{1.3.20}
$$
where the second group is generated by the reflection $s_a$. 

The quadrangle $\M_{2,I}=B_1ABA_2$ has 
$$
P(\M_{2,I})=\{\delta_1=\xi_1+\xi_3,\  
\delta_2=\xi_1-\xi_3,\  \delta_3=\xi_2-\xi_3,\ \delta_4=\xi_2+\xi_3 \} 
\tag{1.3.21}
$$
with the Gram matrix $A_{2,I}$ and the lattice Weyl vector 
$$
\rho_{2,I}=(1/2)\xi_1+(1/2)\xi_2.
\tag{1.3.22}
$$ 
The lattice $M_{2,I}=[a,2b,c]$ 
is 
$$
M_{2,I}\cong 
\langle -8 \rangle \oplus \langle 2 \rangle \oplus \langle 2 \rangle. 
\tag{1.3.23} 
$$  
We have  
$$
W_{2,I}=W^{(2)}(M_{2,I}),\ \ \Sym(\M_{2,I})=D_2 
\tag{1.3.24}
$$ 
where the second group 
is generated by the reflections $s_a$ and $s_{2b}$. 

The right quadrangle $\M_{2,II}=AA_2A_3A_4$ has 
$$
P(\M_{2,II})=\{\delta_1=\xi_3,\ 
\delta_2=2\xi_1-\xi_3,\  \delta_3=2\xi_1+2\xi_2-3\xi_3,\  
\delta_4=2\xi_2-\xi_3\}
\tag{1.3.25}
$$
with the Gram matrix $A_{2,II}$ and the lattice Weyl vector 
$$
\rho_{2,II}=(1/2)\xi_1+(1/2)\xi_2-(1/2)\xi_3.
\tag{1.3.26}
$$
The lattice 
$M_{2,II}=[2\xi_1, 2\xi_2, \xi_3]=[2a,b,2c]$ is 
$$
M_{2,II}\cong U(8)\oplus \langle 2 \rangle.
\tag{1.3.27}
$$ 
We have  
$$
W_{2,II}=W^{(2)}(M_{2,II}),\ \ \ \Sym(P(\M_{2,II}))=D_4
\tag{1.3.28}
$$
where the second group is generated by the reflections $s_{2c}$ and 
$s_{2a}$.

The $8$-gon $\M_{2,III}=B_1A_2B_2A_3B_3A_4B_4A$ has 
$$
\split 
P(\M_{2,III})=\{\delta_1=&\xi_2+\xi_3, \delta_2=3\xi_2-\xi_3,   
\delta_3=\xi_1+4\xi_2-3\xi_3, \delta_4=3\xi_1+4\xi_2-5\xi_3,\\ 
\delta_5=&4\xi_1+3\xi_2-5\xi_3, \delta_6=4\xi_1+\xi_2-3\xi_3, 
\delta_7=3\xi_1-\xi_3, \delta_8=\xi_1+\xi_3\}
\endsplit 
\tag{1.3.29}
$$
with the Gram matrix $A_{2,III}$ and the lattice Weyl vector 
$$
\rho_{2,III}=(1/4)\xi_1+(1/4)\xi_2-(1/4)\xi_3.
\tag{1.3.30}
$$
The lattice $M_{2,III}=[a,2b+c,2c]$ 
is  
$$
M_{2,III}\cong \langle -32 \rangle \oplus \langle 2 \rangle \oplus 
\langle 2 \rangle.
\tag{1.3.31}
$$
We have 
$$
W_{2,III}=W^{(2)}(M_{2,III}),\ \  \Sym(P(M_{2,III}))=D_4
\tag{1.3.32}
$$
where the second group is generated by the reflections 
$s_a$ and $s_{2c}$.  
 
All lattices $M_{2,0}$, $M_{2,I}$, $M_{2,II}$, 
$M_{2,III}$ are characterized by their discriminant quadratic forms.

\subsubhead
1.3.3. Generalized Cartan matrices $A_{3,j}$
\endsubsubhead

On Figure 3 $ABC$ is a triangle with angles 
$0$, $\pi/6$ and $\pi/2$ respectively. Considering the orbit of 
$ABC$ by the group generated by the reflection in $BC$, 
we get a triangle $ABA_2$ with the angles $0$, $\pi/3$ and 
$0$ respectively. Considering the orbit of the triangle 
$ABC$ by the group (of type $D_2$) generated by reflections in $BC$ and 
$AC$, we get a quadrangle $ABA_2B_1$ with 
the angles $0$, $\pi/3$, $0$ and $\pi/3$ respectively. 
Considering the orbit of 
the triangle $ABC$ by the group (of type $D_6$) generated by reflections 
in $BC$ and $BA$, we get a right $6$-gon $AA_2A_3A_4A_5A_6$ 
(here $AA_2$ contains $C$) with zero angles and the center $B$. 
Considering images of the center $B$ by reflections in sides of the 
$6$-gon, we get points $B_1$, $B_2$, $B_3$, $B_4$, $B_5$, $B_6$  
which together with vertices of the $6$-gon 
define an $12$-gon $AB_1A_2B_2A_3B_3A_4B_4A_5B_5A_6B_6$ 
with the angles $0$, $\pi/3$, 
$0$, $\pi/3$, $0$, $\pi/3$,  $0$, $\pi/3$, $0$, $\pi/3$, 
$0$, $\pi/3$  respectively. This 
$12$-gon also is the orbit of the triangle $ABB_1$ with the angles 
$0$, $\pi/6$, $\pi/6$ 
respectively by the group (of type $D_6$) generated by reflections 
in his sides $BA$ and $BB_1$. 

The matrix $A_{3,0}$ is the Gram matrix of the triangle 
$\M_{3,0}=ABA_2$. 
The matrix $A_{3,I}$ is the Gram matrix of the quadrangle 
$\M_{3,I}=ABA_2B_1$.   
The matrix $A_{3,II}$ is the Gram matrix of the right $6$-gon   
$\M_{3,II}=AA_2A_3A_4A_5A_6$, and 
$A_{3,III}$ is the Gram matrix 
of the $12$-gon $\M_{3,III}=AB_1A_2B_2A_3B_3A_4B_4A_5B_5A_6B_6$.   
  
\midinsert
$$\vbox{\centerline{\epsfxsize=5in\epsfbox{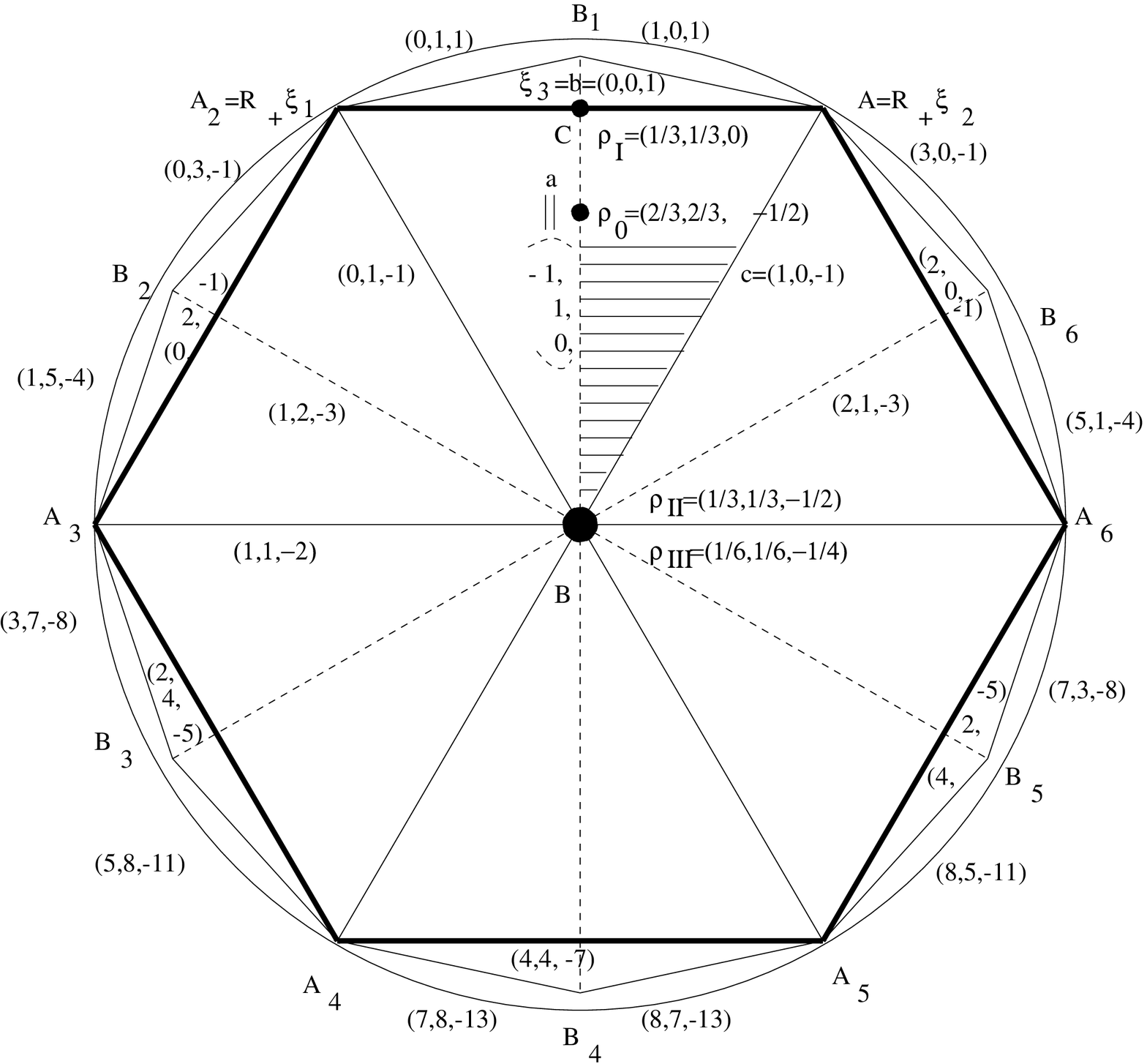}}
\vskip20pt 
\centerline{Figure 3. Fundamental polygons of $A_{3,j}$}}$$
\endinsert

 We consider a hyperbolic lattice 
$M_{3,0}\cong U(3)\oplus \langle 2 \rangle $ with the  
standard bases    
$\xi_1$, $\xi_2$, $\xi_3$ having the Gram matrix 
$$
\left(
\matrix
0 &-3&0\\
-3&0&0\\
0&0& 2
\endmatrix\right). 
$$
We identify the hyperbolic plane of the lattice $M_{3,0}$ and 
of Figure 3:  
Points $A_2=\br_{++}\xi_1$, $A=\br_{++}\xi_2$, elements  
$a=-\xi_1+\xi_2$, $b=\xi_3$ and $c=\xi_1-\xi_3$ are orthogonal 
to sides $BC$, $CA$ and $AB$ respectively of the triangle $ABC$ 
and directed outward. Here $(a,a)=6$ and $(b,b)=(c,c)=2$. 

For $\M_{3,0}=AA_2B$, we have 
$$
P(\M_{3,0})=\{\delta_1=\xi_3,\  
\delta_2=\xi_2-\xi_3,\ \delta_3=\xi_1-\xi_3\}  
\tag{1.3.33}
$$
with the Gram matrix $A_{3,0}$ and the lattice Weyl vector 
$$
\rho_{3,0}=(2/3)\xi_1+(2/3)\xi_2-(1/2)\xi_3.
\tag{1.3.34}
$$    
The lattice $M_{3,0}$, generated by $P(\M_{3,0})$, is 
then $M_{3,0}=[a,b,c]$, and   
$$
M_{3,0}\cong U(3)\oplus \langle 2 \rangle.  
\tag{1.3.35}
$$  
We have 
$$
W_{3,0}=W^{(2)}(M_{3,0}), \ \ \ \Sym(P(\M_{3,0}))=D_1    
\tag{1.3.36}
$$
where the second group is generated by the reflection $s_a$. 

The quadrangle $\M_{3,I}=B_1ABA_2$ has 
$$
P(\M_{3,I})=\{\delta_1=\xi_1+\xi_3,\  
\delta_2=\xi_1-\xi_3,\  \delta_3=\xi_2-\xi_3,\ \delta_4=\xi_2+\xi_3 \} 
\tag{1.3.37}
$$
with the Gram matrix $A_{3,I}$ and the lattice Weyl vector 
$$
\rho_{3,I}=(1/3)\xi_1+(1/3)\xi_2.
\tag{1.3.38}
$$ 
The lattice $M_{3,I}=[a,2b,c]$, 
and it has the discriminant quadratic form 
$$
q_{M_{3,I}}\cong q_1^{(3)}(3)\oplus q_{-1}^{(3)}(3)\oplus q_5^{(2)}(8). 
\tag{1.3.39} 
$$  
We have  
$$
W_{3,I}=W^{(2)}(M_{3,I}),\ \ \Sym(\M_{3,I})=D_2 
\tag{1.3.40}
$$ 
where the second group is generated by the 
reflections $s_a$ and $s_{2b}$. 

The right $6$-gon $\M_{3,II}=AA_2A_3A_4A_5A_6$ has 
$$
\split 
P(\M_{3,II})=\{\delta_1=&\xi_3,\ 
\delta_2=2\xi_2-\xi_3,\  \delta_3=2\xi_1+4\xi_2-5\xi_3,\\  
\delta_4=&4\xi_1+4\xi_2-7\xi_3, 
\delta_5=4\xi_1+2\xi_2-5\xi_3, 
\delta_6=2\xi_1-\xi_3\}
\endsplit 
\tag{1.3.41}
$$
with the Gram matrix $A_{3,II}$ and the lattice Weyl vector 
$$
\rho_{3,II}=(1/3)\xi_1+(1/3)\xi_2-(1/2)\xi_3.
\tag{1.3.42}
$$
The lattice 
$M_{3,II}=[2\xi_1, 2\xi_2, \xi_3]=[2a,b,2c]$ is 
$$
M_{3,II}\cong U(12)\oplus \langle 2 \rangle.
\tag{1.3.43}
$$ 
We have  
$$
W_{3,II}=W^{(2)}(M_{3,II}),\ \ \ \Sym(P(\M_{3,II}))=D_6
\tag{1.3.44}
$$
where the second group is generated by the reflections $s_{2c}$ and 
$s_{2a}$.

The $12$-gon $\M_{3,III}=B_1A_2B_2A_3B_3A_4B_4A_5B_5A_6B_6A$ has 
$$
\split 
P(\M_{3,III})=\{\delta_1=&\xi_2+\xi_3, 
\delta_2=3\xi_2-\xi_3,   
\delta_3=\xi_1+5\xi_2-4\xi_3,\\ 
\delta_4=&3\xi_1+7\xi_2-8\xi_3, 
\delta_5=5\xi_1+8\xi_2-11\xi_3, 
\delta_6=7\xi_1+8\xi_2-13\xi_3,\\ 
\delta_7=&8\xi_1+7\xi_2-13\xi_3, 
\delta_8=8\xi_1+5\xi_2-11\xi_3, 
\delta_9=7\xi_1+3\xi_3-8\xi_3,\\ 
\delta_{10}=&5\xi_1+\xi_2-4\xi_3, 
\delta_{11}=3\xi_1-\xi_3, 
\delta_{12}=\xi_1+\xi_3\}
\endsplit 
\tag{1.3.45}
$$
with the Gram matrix $A_{3,III}$ and the lattice Weyl vector 
$$
\rho_{3,III}=(1/6)\xi_1+(1/6)\xi_2-(1/4)\xi_3.
\tag{1.3.46}
$$
The lattice $M_{3,III}=M_{3,I}=[a,2b,c]$ has the discriminant quadratic form  
$$
q_{M_{3,III}}\cong q_1^{(3)}(3)\oplus q_{-1}^{(3)}(3)\oplus q_5^{(2)}(8). 
\tag{1.3.47} 
$$ 
We have 
$$
W_{3,III} \triangleleft W^{(2)}(M_{3,III})\ \text{and}\   
W^{(2)}(M_{3,III})/W_{3,III}=D_3 
\tag{1.3.48}
$$
where the second group 
is generated by the reflections $s_{\xi_2-\xi_3}$ (in $BA_2$) 
and $s_c$ (in $BA$).  
Moreover, $W_{3,III}$ is normal in a bigger subgroup 
$$
W_{3,III}\triangleleft W^{(2,6)}(M_{3,III})\ \text{and}\ 
W^{(2,6)}(M_{3,III})/W_{3,III}=D_6
\tag{1.3.49}
$$  
where the group $W^{(2,6)}(M_{3,III})$  
is generated by the reflections $s_a$, $s_c$ and $s_{\xi_1+\xi_3}$ 
in all sides of the triangle $B_1BA$. The group 
$$
A(P(\M_{3,III}))=D_6
\tag{1.3.50}
$$
is generated by the reflections $s_a$ and $s_c$.   

All lattices $M_{3,0}$, $M_{3,I}=M_{3,III}$ and $M_{3,II}$  
are characterized by their discriminant quadratic forms.

\head
2. Reflective automorphic Forms and reflective Lattices, and Arithmetic 
Mirror Symmetry Conjecture. 
Lie reflective Automorphic Forms. General Theory 
\endhead 

In this Section 
we want to generalize some results of our papers \cite{N9}, 
\cite{N10}, \cite{GN1} --- \cite{GN4} and \cite{N11}. The 
results of this section were 
subjects of the lectures given by us at the 
meeting of Moscow Mathematical Society  
and Steklov Mathematical Institute (April 1996), 
European Algebraic Geometry meeting at Warwick (July 1996), 
RIMS at Kyoto and Nagoya University (September 1996).   

The main subject of this section is to define proper analogs of 
elliptic and parabolic reflection  groups and corresponding root 
systems (related with lattices of signature $(n,1)$) for lattices with $2$ 
negative squares. These two notions are related by arithmetic 
mirror symmetry conjecture. This section also 
contains our general point of view  
to the theory of automorphic forms related with Mirror Symmetry and 
the theory of automorphic Lorentzian Kac--Moody algebras.      

\subhead
2.1. Reflective hyperbolic lattices and reflective lattices with $2$ negative 
squa\-res. Reflective automorphic forms on IV type domains  
\endsubhead

Let $K$ be a lattice with form $(\cdot ,\cdot )$.  
An element $\delta \in K$ 
is called a {\it primitive root} if $\delta$ is primitive in $K$, 
$(\delta, \delta )>0$ and the 
reflection 
$s_\delta: x \mapsto x-(2(x,\delta)/(\delta,\delta))\delta$ 
in $\delta$ belongs to the automorphism group $O(K)$ 
of the lattice $K$. Equivalently,  
$(\delta,\delta)\,|\,2(K,\delta )$. 
Remark that any element $\delta \in K$ with square 
$(\delta , \delta )=1$ or $2$ is a primitive root. It is called 
$1$ or $2$ root respectively.  
The squares $(\delta, \delta)$ of primitive roots are bounded: 
$(\delta, \delta)\le N$ where $N$ is some constant depending on the 
lattice $K$.  We denote 
by $\Delta (K)_{pr}$ the set of all primitive roots of a lattice $K$.   
A subset $\Delta \subset \Delta (K)_{pr}$ is 
called a {\it primitive root system} if $\Delta $ is 
invariant with respect to the   
reflection group $W(\Delta )$ generated by all reflections in $\Delta$.  
If $\Delta \subset K$ is a root system and 
$K^\prime \subset K$ is a sublattice, then 
$\Delta \cap K^\prime$ is a root system in $K^\prime$. 

\smallpagebreak 

First, we recall some definitions and results 
about reflective hyperbolic lattices and 
corresponding reflection groups and root systems (see \cite{N10}, 
\cite{N4}, \cite{N5}, \cite{V2}, \cite{V3}).  
If $M$ is a hyperbolic lattice (i.e. $M$ has signature $(n,1)$), 
a primitive root system 
$\Delta\subset \Delta (M)_{pr}$ of $M$ is equivalent 
to a reflection group $W(\Delta)\subset O^+(M)$ 
generated by reflections in $\Delta$: 
An element $w\in W(\Delta)$ is a reflection with respect 
to some primitive element $\delta \in M$ if and only if $\delta \in \Delta$. 
A reflection group $W$ of a hyperbolic lattice $M$ 
is called {\it elliptic} if index $[O(M):W]$ is finite. This 
is equivalent to finiteness of volume of the fundamental polyhedron $\M$ 
of $W$ or to the restricted arithmetic type of $W$ together with 
existence of a generalized lattice Weyl vector $\rho$ with square 
$(\rho, \rho)<0$ for the set $P(\M)$ of orthogonal 
vectors to $\M$ (see Definition 1.1.3). 
A reflection group $W$ of a hyperbolic lattice $M$ is called 
{\it parabolic} if $W$ has restricted arithmetic type and the 
set $P(\M)_{pr}$ of primitive orthogonal vectors to $\M$ has a 
generalized lattice Weyl vector $\rho$ with zero 
square $(\rho, \rho)=0$ and does not have a generalized lattice Weyl vector 
with negative square. Using this, we say that a primitive root system 
$\Delta\subset \Delta (M)_{pr}$ is {\it elliptic} or {\it parabolic} 
if the corresponding reflection group $W(\Delta)$ is elliptic 
or parabolic respectively. 
A lattice $M$ is called {\it reflective} if $O(M)$ has at least one 
elliptic or parabolic reflection group $W\subset O^+(M)$ or the 
corresponding elliptic or parabolic primitive root system. 
If the group $W$ is elliptic or parabolic, 
one can also call $M$ {\it elliptic} or {\it parabolic reflective} 
respectively. The main result we need is 
(see \cite{N10}, \cite{N4}, \cite{N5} and 
\cite{V2}) 

\proclaim{Theorem 2.1.1} For $\rk M \ge 3$ the number of reflective 
hyperbolic lattices (elliptic and parabolic) is finite up to 
multiplication of the form of a lattice on $n \in \bq$.  
In particular, the number of maximal 
elliptic and parabolic primitive root systems 
$\Delta$ and corresponding elliptic and parabolic 
reflection groups  $W$ is finite.
\endproclaim  

\smallpagebreak 

We say that a lattice $T$ {\it has $2$ negative squares} if  
$T$ has signature $(n,2)$. We want to suggest a definition of 
reflective lattices with $2$ negative squares and formulate a  
conjectured analog of Theorem 2.1.1. 

Let $T$ be a lattice with $2$ negative squares. 
We denote by $\Omega^+(T)$
a Hermitian symmetric domain of type IV
which is one of the two connected components of
$$
\Omega (T)=\{ Z \in {\Bbb P}(T\otimes \bc) \, \mid \,
(Z, Z)=0,\ (Z, \overline{Z})<0\}.
$$
We consider the corresponding
homogeneous cones (without zeros)
$\widetilde{\Omega} (T)\subset T\otimes \bc$
and
$\widetilde{\Omega}^+(T) \subset T\otimes \bc$ such that
$\Omega (T)=\widetilde{\Omega} (T)/\bc^\ast$ and
$\Omega^+(T)=\widetilde{\Omega}^+(T)/\bc^\ast$.
Let $O^+(T)$ be the subgroup of index two of $O(T)$ which
keeps the component $\Omega^+ (T)$. It is wellknown that
$O^+(T)$ is discrete in $\Omega^+(T)$ and has a fundamental
domain of finite volume.

A function $\Phi$ on $\widetilde{\Omega}^+(T)$ is called
an {\it automorphic form of weight $k$} ($k \in \bz/2$) if $\Phi$ is
holomorphic on $\widetilde{\Omega}^+(T)$,
$\Phi (c\omega)=c^{-k}\Phi (\omega)$ for any $c\in \bc^\ast$,
$\omega \in \widetilde{\Omega}^+(T)$,
and $\Phi (\gamma (\omega))$=
$\chi (\gamma ) \Phi (\omega)$ for any $\gamma \in G$,
$\omega \in \widetilde{\Omega}^+(T)$.
Here $G\subset O(T)^+$ is a subgroup of finite index and
$\chi:G\to \bc^\ast$ is some character or a multiplier system 
with the kernel of finite
index in $G$. Then $\Phi$ is called {\it automorphic with
respect to $G$ with $\chi$}.
The function $\Phi$ additionally should be holomorphic at 
infinity $\Omega (T)_\infty$ of $\Omega (T)$. 
If $\text{codim}_{\Omega (T)}\Omega(T)_\infty\ge 2$, this condition is
automatically valid according to the Koecher principle. Since 
$\dim \Omega (T)_\infty \le 1$, this is in particular true if $\rk T\ge 5$. 

For $e \in T$ with $(e,e) > 0$ we denote
$\Ha_e=\{ Z \in \Omega^+(T)\ \mid \
(Z, e)=0\}$. The $\Ha_e$ is called a {\it
quadratic divisor orthogonal to $e$}. The quadratic
divisor $\Ha_e$ does not change
if one changes $e$ to $te$, $t \in \bq$. The union 
$D(T)$ of all quadratic divisors orthogonal to the primitive roots 
$\Delta(T)_{pr}$ is called {\it discriminant} of $T$. 

\definition{Definition 2.1.2} Let $T$ be a lattice with $2$ negative 
squares. An automorphic form $\Phi$ on $\Omega (T)$ is called 
{\it reflective} for $T$ if the set of zeros 
$(\Phi)_0$ of $\Phi$ is the union 
of quadratic divisors orthogonal to primitive roots 
$\delta \in \Delta (T)$ (with some multiplicities). 
The lattice $T$ is called {\it reflective} if 
there exists at least one reflective automorphic form for $T$. 

We denote 
by $\Delta (\Phi)$ the set of all primitive roots $\delta \in \Delta (T)$ 
such that $\Phi$ is equal to zero on $\Ha_\delta$. Evidently, if 
$\delta \in \Delta (\Phi )$, then $-\delta \in \Delta (\Phi)$. 
An automorphic form $\Phi$ is called 
{\it strongly reflective} if $\Delta (\Phi)$ is a root system and 
$\Phi$ is automorphic with respect to the group $W(\Delta (\Phi))$ 
generated by reflections in all elements of $\Delta (\Phi)$. 

The union 
${D}(\Phi)$ of all quadratic divisors orthogonal to 
$\delta \in \Delta(\Phi)$ is called {\it discriminant} of $\Phi$. 
It is why reflective automorphic forms are also called 
{\it discriminant automorphic forms}. They are related  
with discriminants of K3 surfaces moduli (see \cite{GN3}, \cite{N11}).

Below we will also use notation ${D}(\Delta)$ for the union of 
all quadratic divisors orthogonal to elements of some subset 
$\Delta\subset \Delta (T)_{pr}$ of primitive roots.     
\enddefinition

Suppose that $\Phi$ is a reflective automorphic form 
for $T$. Suppose that $\Phi$ 
is invariant with respect to a 
subgroup $G\subset O^+(T)$ of a finite index. The product  
$$
\Phi^\prime = \prod_{g \in G\backslash O^+(T)}{g^\ast \Phi }, 
$$
is also a reflective automorphic form, but $\Delta(\Phi^\prime)$ 
is invariant with respect to $O^+(T)$. 
In particular, $\Delta (\Phi^\prime)$ is a root system. 
Thus, $\Phi^\prime$ is a strongly reflective automorphic form. 
Thus, studying reflective lattices, one can restrict 
himself considering strongly reflective automorphic forms. 
The advantage of strongly reflective automorphic forms 
is that they are automorphic with some quadratic character with 
respect to the reflection group $W(\Delta (\Phi))$ generated by 
the reflections in $\Delta (\Phi)$.  
We remark that $W(\Delta (\Phi))$ always has finite index in 
$O(T)$ because of the general results by G.A. Margulis about 
normal subgroups of arithmetic groups of rank $\ge 2$.

It is a very interesting and difficult problem to find all 
reflective lattices $T$ with $2$ negative squares and all their  
reflective automorphic forms.  
Some reflective 
lattices $T$ and their reflective automorphic forms 
were found in \cite{B1}---\cite{B6} and 
\cite{GN1}-- \cite{GN4}. Many other will be found in Part II of 
this paper.  

\subhead
2.2. Arithmetic Mirror Symmetry Conjectures  
\endsubhead 

Conjecture below was first formulated in \cite{N11}  
for 2-roots and corresponding $2$-reflective automorphic forms. 
We consider this conjecture as a {\it mirror symmetric statement} 
to Theorem 2.1.1. 

\proclaim{Conjecture 2.2.1} The number of reflective 
lattices $T$ with $2$ negative squares and $\rk T\ge 5$  
is finite up to multiplication of the form of the lattice $T$. In 
particular, for $\rk T\ge 5$ the number of maximal systems of 
primitive roots $\Delta (\Phi)$ 
or reflective automorphic forms $\Phi$ is finite up to multiplication of 
the form of $T$.   
\endproclaim 

The main reason to have so strong statement is the Koecher principle 
for automorphic forms (see \cite{Ba}, for example). 
 
Let $T$ be a reflective lattice with $2$ negative squares and 
$\Phi$ be a reflective automorphic form for $T$. 
Let $\Delta (\Phi )$ be the 
set of primitive roots of $\Phi$. 
Let $T_1\subset T$ be a primitive sublattice of $T$ with $2$ 
negative squares.  
Then $\Omega (T_1)\subset \Omega (T)$ is a symmetric subdomain of type 
IV of $\Omega (T)$ and restriction $\Phi |_{\Omega (T_1)}$ 
is an automorphic form 
for the lattice $T_1$. By Koecher principle, 
$\Phi|_{\Omega (T_1)}$ should have a non-empty set of zeros if 
$\text{codim}_{\Omega(T_1)}\Omega (T_1)_\infty \ge 2$. Thus  
$$
{D}(\Delta(\Phi))\cap \Omega (T_1)\not=\emptyset\ 
\text{if}\  \text{codim}_{\Omega(T_1)}\Omega (T_1)_\infty \ge 2 .
$$
We get the following condition:  

\proclaim{Necessary Condition 2.2.2} Let $T$ be a lattice with $2$ 
negative squares and $\Delta \subset \Delta (T)_{pr}$ a set of primitive 
roots. 
If there exists a reflective automorphic form $\Phi$ with   
the set of roots $\Delta (\Phi)\subset \Delta$, then 
for any primitive sublattice 
$T_1\subset T$ with $2$ negative squares and  
$\text{codim}_{\Omega(T_1)}\Omega (T_1)_\infty \ge 2$,  
one has 
$$
{D}(\Delta )\cap \Omega (T_1)\not=\emptyset. 
$$

In particular, for $\Delta=\Delta(T)_{pr}$ we get: 
If $T$ is a reflective lattice, then  
for any primitive sublattice $T_1\subset T$ with $2$ negative 
squares and $\text{codim}_{\Omega(T_1)}\Omega (T_1)_\infty \ge 2$, 
the set  
${D}(\Delta (T)_{pr})\cap \Omega (T_1)$ is not empty 
\endproclaim

It was shown in \cite{N11} that Necessary Condition 2.2.2 for a lattice 
$T$ to be reflective is extremely strong. Actually, we know only 
one way to prove that a lattice $T$ of $\rk T\ge 5$ 
satisfies this condition: to construct a reflective automorphic form 
for $T$. It is why we suggest more strong than Conjecture 2.2.1   

\proclaim{Conjecture 2.2.3} The set of lattices $T$ with $2$ negative 
squares and with $\rk T\ge 5$ 
which satisfy Necessary Condition 2.2.2 is finite up to multiplication of 
the form of $T$ on $n \in \bq$. In particular, the set of maximal 
systems of primitive roots  
$\Delta$ which satisfy Necessary Condition 2.2.2 is finite 
up to multiplication of the form on $n \in \bq$.  
\endproclaim
 
\smallpagebreak 

Now we have two types of reflective lattices: hyperbolic reflective 
lattices and reflective lattices with $2$ negative squares. We think 
that these two types of lattices are related by  
{Arithmetic Mirror Symmetry Conjecture} which we formulate below. 
It is inspired by Mirror Symmetry (at least, for K3 surfaces).  
We consider this conjecture as 
much more strong variant of our Mirror Symmetry 
Conjectures in \cite{GN3} and \cite{N11}. 
But the main ingredient of our Conjecture in 
\cite{GN3} is preserved. We believe that the right variant of mirror symmetry 
for K3 surfaces is connected with the very global objects like 
Automorphic Forms on a IV type domains. Finiteness Conjectures 
2.2.1 and 2.2.3.  together with Theorem 2.1.1 show that there actually exist 
very few cases when this kind of theory takes place. We believe that these  
finiteness results and Conjectures are related with conjectured 
finiteness of families of Calabi--Yau 3-folds and $n$-folds   
(one can try to find some connection in \cite{N12}, \cite{HM1}, 
\cite{HM2}, \cite{Ka1}, \cite{Ka2},  \cite{DVV}).    
 
\proclaim{Arithmetic Mirror Symmetry Conjecture 2.2.4} 

(i) Let $T$ be  
a reflective lattice with $2$ negative squares 
and $\Phi$ be a reflective automorphic form for $T$ 
with a primitive root system $\Delta (\Phi)$. Let $c \in T$ be a primitive 
isotropic element of $T$ and $S=c^\perp_T/[c]$ be the corresponding 
hyperbolic lattice with the root system 
$$
\Delta (\Phi)|_S=\Delta (\Phi)\cap c^\perp_T~\mod [c].
$$
Then the root system $\Delta (\Phi)|_S\subset S$ is elliptic 
or parabolic (in particular, the 
lattice $S$ is reflective) if the set of roots $\Delta (\Phi)|_S$ is 
non-empty. 

(ii) Any hyperbolic reflective lattice $S$ with  
an elliptic or parabolic primitive root system $\Delta$ 
may be obtained from some lattice $T$ with $2$ negative squares 
and a reflective automorphic form 
$\Phi$ of $T$ by the construction (i) above. 
\endproclaim 

We have two arguments to support Conjecture 2.2.4. First we 
will show in Theorem 2.2.5 below that if $T$ together 
with a primitive root system $\Delta$ satisfies the condition 
(i) of Conjecture 2.2.4, then 
$T$ with $\Delta$ satisfies Necessary 
Condition 2.2.2 for ``the most part'' 
of primitive sublattices $T_1\subset T$. We have mentioned 
that this condition is extremely strong which had 
been shown in \cite{N11}.   

\proclaim{Theorem 2.2.5} Let $T$ be a lattice with $2$ negative 
squares and $\Delta$ be a primitive root system of $T$. Suppose 
that for any primitive isotropic $c \in T$ and 
$S=c^\perp_T/[c]$ the root system 
$\Delta |_S\subset S$ is elliptic or parabolic if it is not empty 
(i.e. the condition (i) of Conjecture 2.2.4 is valid for $\Delta$.) 
Then $T$ satisfies Necessary Condition 2.2.3 for all 
primitive sublattices $T_1\subset T$ such that $T_1$ has a primitive 
isotropic element $e$ with property $e^\perp_T\cap \Delta \not=\emptyset$. 
\endproclaim

\demo{Proof} Let $S=e^\perp_T/[e]$ and let $S_1=e^\perp_{T_1}/[e]$ be 
its hyperbolic sublattice. 
It is sufficient to show that the root system 
$\Delta |_S$ has the following property which is a hyperbolic analog of 
Necessary Condition 2.2.2: 
$$
(\bigcup_{\delta \in \Delta |_S}{\Ha_\delta})\cap \La (S_1)\not=\emptyset .
\tag{2.2.1}
$$
Here $\La (S_1)\subset \La (S)$ are the hyperbolic spaces 
corresponding to $S_1\subset S$, and $\Ha_\delta$ is the 
hyperplane in $\La (S)$ orthogonal to $\delta$.  
It is true that $\rk S_1\ge 2$ and $S_1$ does not have 
isotropic elements if $\rk S_1=2$ (it follows from the condition 
$\text{codim}_{\Omega(T_1)}\Omega (T_1)_\infty \ge 2$ 
where we consider an empty set as having dimension $(-1)$). It follows 
that the subspace $\La (S_1)$ has a non-empty finite (not at infinity) 
intersection with a face of one of fundamental polyhedron $\M$ for 
$W(\Delta |_S)$ acting in $\La (S)$. 
This is obvious if $\M$ has finite volume (elliptic case). 
One should draw a picture 
to see that the same is true for a parabolic case when the polyhedron 
$\M$ is parabolic with respect to a cusp in $S$. 
This finishes the proof.  
\enddemo

Another support of Conjecture 2.2.4 is that if a reflective 
automorphic form $\Phi$ has a Fourier expansion with a generalized 
lattice Weyl vector (one can consider this as a special type 
of Fourier expansion), then the statement (i) of Conjecture 2.2.4 
is valid for this cusp (see Sect. 2.4). 
Infinite products ``a la Borcherds'' (see R. Borcherds, \cite{B5}) 
give a particular case of Fourier expansions 
with a generalized lattice Weyl vector and are related with 
automorphic forms with rational quadratic divisors (see \cite{B5}). 
It is possible that reflective automorphic forms $\Phi$ always have 
these infinite product expansions at cusps because of their special 
sets of zeros. One can also believe that 
this is true because of finiteness Conjectures 2.2.1 and 2.2.3. 
It seems, all known reflective automorphic forms $\Phi$ have 
these infinite product expansions at cusps.   

\subhead
2.3. Fourier expansion of an automorphic form $\Phi$ on 
the domain $\Omega(T)$ 
\endsubhead

Let $T$ be  a lattice with $2$ negative squares.
Here we analyze Fourier  expansion of 
an automorphic form $\Phi$ on $\Omega^+(T)$ at a 
$0$-dimensional  cusp 
in terms of the  lattice $T$. We describe canonical linear complex
coordinates at the cusp and the group of translations
which determines the Fourier expansion.

A cusp (of dimension $0$) of $\Omega^+(T)$ is defined by an isotropic 
sublattice of $T$ of rank one, i.e. by a 
primitive isotropic element $c \in T$
which  is defined by the cusp canonically up to replacing $c$ to $-c$.
For the  given isotropic primitive $c \in T$ 
we put  
$$
S_c= c^\perp/ c:=c^\perp_T/\bz c
$$
which is  a  hyperbolic lattice canonically defined by the cusp.
Let 
$O^+(T)_c=\{g \in O^+(T) \ |\ g(c) = c\}$
be the corresponding  to $c$ parabolic subgroup
and let
$$
O^+(T)_{c^\perp/c }=\{\phi \in O^+(T)_c\ |\ \phi \
\text{is identical on } c^\perp/c \}
$$
be its unipotent subgroup which is a subgroup of finite index 
of the additive group $S_c^\ast\subset S\otimes \bq$. 
The corresponding embedding $f$
is defined as follows:
for $\phi \in O^+(T)_{c^\perp /c}$ the function
$f(\phi)\in S_c^\ast$ is determined by the rule:
$$
\phi (x)=x+{f(\phi)}(x)c, \ \ x \in c^\perp.
$$
One can prove that $f(\phi)$ can be 
extended to a function on $T$ from the module
$$
T^\ast_{c,0}=\{x \in T^\ast\ |\ (x,c)=0\}.
$$
It means that in fact
$$
O^+(T)_{c^\perp/c }\subset \widetilde{T}^\ast_{c,0}
\subset S_c^\ast
\tag{2.3.1}
$$
where  $T^\ast_{c,0}\to \widetilde{T}^\ast_{c,0}
\subset S_c^\ast$ is the natural homomorphism.

Let us define canonical coordinates at the cusp $c$.
For an arbitrary element $\bc \omega \in \Omega (T)$ one can choose
canonically its representative $\omega \in \bc \omega$ by the
condition $(\omega, c)=-1$. This choice  identifies
$\Omega (T)$ with a domain of the affine quadric
$$
\Omega (T)_c=\{ \omega \in T\otimes \bc \, \mid \,
(\omega, \omega)=0,\ (\omega, \overline{\omega})<0,\ (\omega,c)=-1\}
\cong \Omega (T)
$$
and $\Omega^+(T)$ with a connected component $\Omega^+(T)_c$ of 
$\Omega (T)_c$. 
We consider  a cone 
$$
V (S_c)=\{x \in S_c\otimes \br \ |\ (x,x)<0\}
$$ 
of the hyperbolic lattice $S_c$,
its half-cone
$V^+(S_c)$ and the corresponding complexified cone
$\Omega (V^+(S_c))=S_c\otimes \br+i V^+(S_c)$
which is a complex tube domain of type IV.
An element  
$$
h\in T^\ast_{c,-1}=\{x \in T^\ast \ |\ (x,c)=-1\} 
$$
defines an orthogonal decomposition
$
T\otimes \bq=S_c\otimes \bq\oplus (\bq c+\bq h) 
$
of $T$ over $\bq$
and an isomorphism (after 
replacing the cone $V^+(S_c)$ by the opposite half-cone 
$-V^+(S_c)$ if necessary)
$$
\Omega (V^+(S_c))\to \Omega^+(T)_c\cong \Omega^+(T)
\ \text{ by }\ 
z \mapsto z\oplus 
(\tsize\frac{(z,z)+(h,h)}{2}c+h)\in \Omega^+(T)_c.
$$
The coordinate $z \in \Omega (V^+(S_c))$
depends on the choice of $h \in T^\ast_{c,-1}$.
If one replaces $h$ with another $h^\prime \in T^\ast_{c,-1}$,
then $h^\prime-h=r\in T^\ast_{c,0}$.
The new coordinate $z^\prime \in \Omega (V^+(S_c))$ is then
$
z^\prime =z-\widetilde{r}
$
where $\tilde {r}$ is the natural image of $r$ in $S_c^\ast$
(see \thetag{2.3.1}). 

Using the coordinate $z \in \Omega (V^+(S_c))$ defined 
by the choice of $h \in T^\ast_{c,-1}$, we can identify
an arbitrary automorphic form $\Phi$ of  weight $k$
on the domain $\Omega^+(T)$
with an {\it automorphic form $\Phi_{c,h}$ on the
domain $\Omega (V^+(S_c))$}. 
They are related by the formula:
$$
\Phi \bigl(\lambda(z\oplus (\tsize\frac{(z,z)+(h,h)}2c+h))\bigr)=
\lambda^{-k}\Phi_{c,h}(z),\quad z \in \Omega (V^+(S_c)),\ 
\lambda \in \bc^\ast. 
$$
In other words, after the $h$-identification of the domains
$\Omega (V^+(S_c))=\Omega^+(T)_c$,  
we get  $\Phi_{c,h}=\Phi |_{ \Omega^+(T)_c}$.
For the  coordinates defined by  $h, h^\prime \in T^\ast_{c,-1}$ 
with  $r=h^\prime-h \in T^\ast_{c,0}$,   
one has $\Phi_{c,h^\prime}(z^\prime)=\Phi_{c,h}(z)$   
if  $z^\prime  =z-\widetilde{r}$. Thus,
$$
\Phi_{c,h^\prime }(z)=\Phi_{c,h} (z+\widetilde{r})
\tag{2.3.2}
$$
for any  $z \in \Omega (V^+(S_c))$, any  $h,h^\prime \in T^\ast_{c,-1}$
with $r=h^\prime-h$
and 
$\widetilde{r}\in \widetilde{T}^\ast_{c,0} \subset S_c^\ast$.

Let us assume that $\Phi$ is an automorphic form with respect
to a group $G\subset O^+(T)$ of finite index with a 
multiplier system $\chi$.
The unipotent subgroup
$$
G_{c^\perp/c }=\{g \in O^+(T)_{c^\perp/c}\cap G\ |\ 
g^*\Phi=\Phi \}
$$
has   a finite index in the lattice $S_c^\ast$.
The invariance of $\Phi$ with 
respect to the subgroup $G_{c^\perp/c}$ is equivalent to
the periodicity of $\Phi_{c,h}$
$$
\Phi_{c, h}(z+ \alpha)=\Phi_{c,h}(z),\qquad z \in \Omega (V^+(S_c)),\
\alpha \in G_{c^\perp/c}\subset
\widetilde{T}^\ast_{c,0}\subset S_c^\ast.
$$
Thus, one can consider 
the Fourier expansion of $\Phi_{c,h}$
by characters of  $(G_{c^\perp/c})^\ast$.
We may  consider $(G_{c^\perp/c})^\ast$ as a subset
of $S_c^\ast \otimes \bq=S_c\otimes \bq$. Then
$S_c\subset (G_{c^\perp/c})^\ast$.
As a result  we get 
$$
\Phi_{c,h}(z)=\sum_{a \in (G_{c^\perp/c})^\ast}
{m(a)\exp{(-2\pi i (a ,z))}}.
\tag{2.3.3}
$$
This is called {\it the Fourier expansion of the automorphic
form $\Phi$ at the cusp $c$}. In general this Fourier
expansion depends on the  choice of $h \in T_{c,-1}^\ast$ and
according to \thetag{2.3.2} there  is  in fact 
a homogeneous  space over the finite Abelian group
$$
\widetilde{T}^\ast_{c,0}/\,(G_{c^\perp/c})
\tag{2.3.4} 
$$
of Fourier expansions of $\Phi$ for different canonical coordinates. 
All such  Fourier expansions are related according to \thetag{2.3.2}.
For another $h^\prime \in T^\ast_{c,-1}$, 
we get 
$$
\Phi_{c,h^\prime}(z)=\sum_{a \in (G_{c^\perp/c})^\ast}
{\exp{(-2\pi i (a, \widetilde{r}))}\, m(a)
\exp{(-2\pi i (a ,z)})}.
\tag{2.3.5}
$$
where $r=h^\prime -h$, $\widetilde{r}\in \widetilde{T}^\ast_{c,0}$ and
$\exp{(-2\pi i (a, \widetilde{r}))}$ is an element  
from the group of unities 
of the order equals to the exponent of the finite 
Abelian group $\widetilde{T}^\ast_{c,0}/\,(G_{c^\perp/c})$.
Below we will be especially interesting in Fourier expansions with
integral Fourier coefficients $m(a)$. Using \thetag{2.3.5}, one
can find all Fourier expansions (in the different $h$-coordinates) 
of the automorphic form $\Phi$ with integral Fourier coefficients 
from one of them.

{\it Further, when we speak about some properties
of  Fourier expansion  of an automorphic form $\Phi$ at a cusp $c$, 
we mean existence of one of the Fourier expansions of $\Phi$  
at the cusp $c$ with this property.}
Below we  omit the index $h\in T_{c,-1}^\ast$ in 
Fourier expansion of an automorphic form $\Phi$.   

\subhead
2.4. Fourier expansions of reflective automorphic forms and Arithmetic 
Mirror Symmetry Conjecture 
\endsubhead
 
Let $T$ be a lattice with $2$ negative squares and $\Phi $ be an  
automorphic form on $\Omega ^+(T)$ possibly with some character 
with respect to a subgroup $G\subset O^+(T)$ of finite index. 
Let $c\in T$ be a primitive isotropic element and $S=c^\perp/[c]$ 
the corresponding hyperbolic lattice. Then the automorphic form  
$\Phi(z)$, 
$z \in \Omega (V^+(S))$, is invariant with some character $\epsilon$ 
with respect to 
a subgroup $H\subset O^+(S)$ of a finite index. 
Let us suppose that $H$ is a semi-direct product 
$$
H=W\rtimes A
$$
where $W$ is a reflection group with a fundamental polyhedron $\M$ in 
$\La ^+(S)$ and $A$ is a group of symmetries of $\M$. Then $\Phi (z)$ 
has a Fourier expansion  
$$
\Phi (z)=\sum_{w \in W}{\epsilon (w)}\sum_{a \in \br_{++}\M}
m(a)\exp(-2\pi i (w(a), z)).
\tag{2.4.1}
$$
Here $a$ runs through some lattice $M$ in $S\otimes \bq$, and 
we call elements $a$ of this lattice {\it Fourier exponents}. The 
corresponding $m(a )\in \bc$ is called {\it Fourier coefficient}.

\definition{Definition 2.4.1} We say that a Fourier exponent  
$\rho \in \br_{++}\M$ is a {\it generalized lattice Weyl vector} for 
$\Phi$ (with respect to $W$ and $\M$) 
if there exists a partial ordering $\le $ of Fourier exponents  
$a \in \br_{++}\M$ such that this ordering is $A$-invariant and  
$\rho$ is the unique minimal element for this ordering.

We remark that $(\rho ,\rho) \le 0$ since $\br_{++}\M \subset 
\overline{V^+(S)}$. 

\enddefinition

\proclaim{Proposition 2.4.2} If there exists a generalized lattice 
Weyl vector $\rho$ for $\Phi$, then the reflection group 
$W$ is elliptic if $(\rho, \rho)<0$, and $W$ is elliptic or parabolic 
if $(\rho, \rho)=0$. In particular, the lattice $S$ is reflective. 
\endproclaim

\demo{Proof} Obviously, $\rho$ is fixed by $A$. If $(\rho, \rho)<0$, then 
$A$ is finite since $A$ is discrete in $\La ^+(S)$. It follows that 
$W$ has finite index in $O^+(S)$ because $W\rtimes A$ does. Thus, $W$ 
is elliptic. 
If $(\rho, \rho)=0$, but $A$ is finite, we get the previous case. 
If $A$ is infinite, then $W$ is parabolic with respect to the cusp $\rho$,  
by definition. It follows the statement.  
\enddemo

As an example of a possible ordering for exponents in $\br_{++}\M$, 
we can consider a {\it standard ordering} defined by the cone 
$\overline{V^+(S)}$. Thus, $x\ge y$ if there exists $z \in \overline{V^+(S)}$ 
such that $x=y+z$. For the standard ordering, the Fourier expansion of 
$\Phi (z)$ with a generalized 
lattice Weyl vector $\rho \in \br_{++}\M$ has the form      
$$
\Phi (z)=\sum_{w \in W}{\epsilon (w)}\sum\Sb \rho+a \in \br_{++}\M,\\
a \in \overline{V^+(S)}\endSb 
m(a)\exp(-2\pi i (w(\rho+a), z)).
\tag{2.4.2} 
$$

Let $\Delta\subset S$ be the set of all primitive roots of $W$. 
A root $\delta \in \Delta$ is called {\it positive} if $(\delta, \M)\le 0$. 
Otherwise, $\delta$ is called negative.  
We denote by $\Delta_+$ the set of all positive roots of $\Delta$.  
The most general ordering related with $W$, we can introduce, 
is the ordering defined by the closed convex cone $\goth T$ below: 
$$
x \ge y\ \ \text{if\ }\  x-y \in {\goth T}=
\overline{\sum_{\delta \in \Delta_+}{\br_+\delta} +V^+(S)}.  
\tag{2.4.3}
$$
Fourier expansion of $\Phi (z)$ 
with a generalized lattice Weyl 
vector $\rho \in \br_{++}\M$ corresponding to 
the ordering \thetag{2.4.3} is     
$$
\Phi (z)=\sum_{w \in W}{\epsilon (w)}\sum\Sb \rho+a \in \br_{++}\M,\\
a \in {\goth T}\endSb 
m(a)\exp(-2\pi i (w(\rho+a), z)).
\tag{2.4.4} 
$$

A particular case of Fourier expansion of $\Phi$ of type 
\thetag{2.4.4} with a generalized lattice Weyl vector is given by the 
{\it infinite product} ``a la Borcherds'' (see \cite{B5}). 
This is a formal infinite product expansion of the form 
$$
\Phi (z)=
C\exp{(-2 \pi i (\rho, z))}
\prod_{\alpha >0}
{(1-u(\alpha)\exp(-2\pi i (\alpha, z))^{\mult~ \alpha}}. 
\tag{2.4.5}
$$
Here $C$ is a constant, $\rho \in \br_{++}\M$, $\alpha >0$ means that 
$\alpha \in \bn \Delta_+\cup (M\cap \overline{V^+(S)})$,   
coefficients $u(\alpha)\in \bc$ and 
multiplicities $\mult~ \alpha\in \bq$. Obviously, then $\Phi (z)$ 
has Fourier expansion of type \thetag{2.4.4} with 
the generalized lattice Weyl vector $\rho\in \br_{++}\M$. 

Thus, from considerations above and Proposition 2.4.2, we get  

\proclaim{Proposition 2.4.3} Let an automorphic form $\Phi$ has 
a Fourier expansion of type \thetag{2.4.4} with a generalized 
lattice Weyl vector 
$\rho \in \br_{++}\M$. Then the reflection group $W$ is  
elliptic or parabolic and the lattice $S$ is reflective.  
In particular, this is valid if $\Phi (z)$ has a formal infinite product 
expansion of type \thetag{2.4.5}.
\endproclaim 

From Proposition 2.4.3, we get the statement which supports 
Arithmetic Mirror Symmetry Conjecture 2.2.4.  

\proclaim{Corollary 2.4.4} Let $T$ be a lattice with $2$ negative 
squares and $\Phi$ a reflective automorphic form with primitive roots 
system $\Delta (\Phi )$. Suppose that at a cusp $c$ the form $\Phi$ 
has a Fourier expansion of type \thetag{2.4.4} with a generalized 
lattice Weyl vector 
$\rho \in \br_{++}\M$ for $W=W((\Delta (\Phi)\cap c^\perp) \mod [c])$ 
(for example, assume that $\Phi$ has a formal infinite product expansion 
\thetag{2.4.5}). 
Then Arithmetic Mirror Symmetry Conjecture 2.2.4 is valid for $\Phi$ 
at the cusp $c$.   
\endproclaim 

It is possible that all reflective automorphic forms have  
infinite product expansions of type \thetag{2.4.5} at cusps  
because of their very special sets of zeros.   
 
\subhead
2.5. Lie reflective automorphic forms 
\endsubhead

In notation above, let $P(\M)\subset S$ be an acceptable set of orthogonal 
vectors to $\M$. It means that 
$$
(\alpha, \alpha )\,|\,2(\alpha, \alpha^\prime)\ \text{for all\ }
\alpha, \alpha^\prime \in P(\M ). 
\tag{2.5.1}
$$
We choose a subset $P(\M)_\1o\subset P(\M)$ which is called  
{\it odd subset}. Its complement $P(\M )_\0o=P(\M)-P(\M)_\1o$ is 
called {\it even subset}. We additionally suppose that 
$$
(\alpha, \alpha)\,|\,(\alpha, \alpha^\prime),\ \text{for all\ }
\alpha \in P(\M)_\1o,\ \alpha^\prime \in P(\M).
\tag{2.5.2}
$$
The set $P(\M)$ together with its subdivision 
$P(\M)= P(\M)_\0o \coprod P(\M)_\1o$ on even and odd subsets 
satisfying \thetag{2.5.1} and \thetag{2.5.2} is called {\it acceptable}. 
We remark that the main invariant of $W$, $\M$ and $P(\M)$ is the 
{\it symmetrized generalized Cartan matrix}
$$
B=\left((\alpha, \alpha^\prime)\right),\ \ \alpha, \alpha^\prime \in P(\M),
\tag{2.5.3} 
$$ 
or (if one permits to multiply the matrix $B$ by constants) 
a generalized Cartan matrix  
$$
A= \left({2(\alpha, \alpha^\prime)\over(\alpha, \alpha)}\right), 
\tag{2.5.4}
$$
and the subset of odd indices $P(\M )_\1o\subset P(\M)$.  

\definition{Definition 2.5.1} We say that $\Phi$ has Fourier 
expansion at the cusp $c$ of {\it Lie type 
(with the symmetrized generalized Cartan matrix $B$ or with the 
generalized Cartan matrix $A$)} if its Fourier expansion has the form  
(it is a very special type of the expansion \thetag{2.4.2} 
which corresponds to the subcone $\br_{+}\M\subset \overline{V^+(S)}$)  
$$
\Phi (z)=\sum_{w \in W}{\epsilon (w)}\left(
\exp(-2\pi i (w(\rho),z))-
\sum\Sb a \in M \cap \br_{++}\M\endSb 
m(a)\exp(-2\pi i (w(\rho+a), z))\right) 
\tag{2.5.5} 
$$
where 

\noindent 
(i) $M\subset S\otimes \bq$ is a lattice in $S\otimes \bq$ such that 
$P(\M)\subset M$,  
$(\alpha, \alpha)\,|\,2(M, \alpha)$ for any $\alpha \in P(\M)$, and 
$(\alpha, \alpha)\,|\,(M, \alpha)$ for any $\alpha \in P(\M)_\1o$; 
the lattice $M$ is called a {\it root lattice}.    

\noindent 
(ii) $\epsilon (s_\alpha )=(-1)^{1+\io}$ for any  
$\alpha \in P(\M)_{\io}$ and any $\io=\0o, \1o$. 

\noindent 
(iii) The element $\rho \in \br_+\M \cap (S\otimes \bq )$ 
satisfies the equality:       
$(\rho, \alpha )=-(\alpha, \alpha)/2$ for any $\alpha \in P(\M)$;   
it is called a {\it lattice Weyl vector}. 
 
\noindent 
(iv) All Fourier coefficients $m(a)$, 
$a \in M\cap \br_{++}\M$, are integral.    
\enddefinition

Using construction \cite{K1}, \cite{K3}, \cite{B1}, \cite{GN1},  
\cite{R}, one can correspond to Lie type Fourier expansion \thetag{2.5.5} 
a generalized Kac--Moody superalgebra. It is defined by a set of 
simple roots ${}_s\Delta$ which is divided on the set of simple 
real roots ${}_s\Delta^\re$ and the set of simple imaginary roots 
${}_s\Delta^\im$. Both these sets are divided in subsets 
of even and odd roots which are marked by indices $\0o$ and $\1o$ 
respectively. We have 
$$
{}_s\Delta^\re_\0o=P(\M)_\0o,\ \ 
{}_s\Delta^\re_\1o =P(\M)_\1o 
\tag{2.5.6} 
$$ 
and   
$$
{}_s\Delta^\im ={}_s\Delta^\im_\0o \cup  
{}_s\Delta^\im_\1o,\ \ {}_s\Delta^\im_\io= 
\{m(a)_\io a \ | a \in M\cap \br_{++}\M\}
\tag{2.5.7}  
$$
where non-negative integers 
$m(a)_\0o$ and $m(a)_\1o$ are related with Fourier coefficients 
$m(a)$ in \thetag{2.5.5} as follows:  
$$
m(a)_\0o\  - m(a)_\1o=m(a),\ \text{for any $a \in \br_{++}\M$ with\ }(a,a)<0,
\tag{2.5.8} 
$$ 
and for any primitive $a_0\in M\cap \br_{++}\M$ with 
$(a_0,a_0)=0$ one has the identity of formal power series with 
one variable $q$:  
$$
\prod_{n \in \bn}{(1-q^n)^{m(na_0)_\0o\ -m(na_0)_\1o}}=1-
\sum_{t \in \bn}{m(ta_0)q^t}.
\tag{2.5.9} 
$$
In  \thetag{2.5.7} $m(a)_\io a$ means that we repeat the element $a$ 
in the sequence ${}_s\Delta^\im_\io$ exactly $m(a)_\io$ times. The 
set ${}_s\Delta={}_s\Delta^\re \cap {}_s\Delta^\im$ is called 
the {\it set of simple roots}.

The {\it generalized Kac--Moody superalgebra} 
$\geg=\geg^{\prime\prime}({}_s\Delta )$ 
corresponding to ${}_s\Delta$ is a Lie superalgebra over $\br$ 
generated by $h_r$, $e_r$, $f_r$, $r \in {}_s\Delta$. All $h_r$ are 
even; $e_r$, $f_r$ are even (respectively odd) if $r$ is even 
(respectively odd).  
The algebra $\geg$ has the defining relations (where $r,r^\prime$ are 
arbitrary elements of ${}_s\Delta$):  

(1) The map $r \mapsto h_r$ gives an embedding
of $M\otimes \br$ into $\geg^{\prime\prime}({}_s\Delta)$ as
an Abelian subalgebra (it is even since all $h_r$ are even).
In particular, all elements $h_r$ commute.

(2) $[h_r, e_{r^\prime}]=(r, r^\prime)e_{r^\prime}$, and
$[h_r, f_{r^\prime}]=-(r,{r^\prime})f_{r^\prime}$.

(3) $[e_r, f_{r^\prime}]=h_r$ if $r=r^\prime$, and is $0$ if
$r \not=r^\prime$.

(4) $(\text{ad\ } e_r)^{1-2(r,r^\prime)/(r,r)}e_{r^\prime }=
(\text{ad\ } f_r)^{1-2(r,r^\prime)/(r,r)}f_{r^\prime }=0\
\text{if $r\not= r^\prime$ and $(r,r)>0$}$ \hfil\hfil
\newline
\phantom{(4)(4)} (equivalently, $r \in {}_s\Delta^{\re}$).

(5) If $(r,r^\prime)=0$, then $[e_r, e_{r^\prime}]=[f_r,f_{r^\prime}]=0$.
 
\vskip5pt

The superalgebra $\geg=\geg ^{\prime\prime}({}_s\Delta)$
is graded by the root lattice $M$ as follows. Let
$$
{Q}_+ =\sum_{\alpha \in {}_s\Delta}{\bz_+\alpha }\subset M
\tag{2.5.10}
$$
be the integral cone (semi-group) generated by all simple roots $\Delta$. 
We have
$$
\geg=\left( \bigoplus_{\alpha \in {Q}_+} {\geg_\alpha} \right)
\bigoplus \geg_0 \bigoplus
\left( \bigoplus_{\alpha \in {Q}_+} {\geg_{-\alpha}} \right)
\tag{2.5.11}
$$
where $e_r$ and $f_r$ have degree $r\in {Q}_+$ and
$-r \in -{Q}_+$ respectively, $r \in {}_s\Delta$;
and $\geg_0=M\otimes \br$.
An element $\alpha \in \pm {Q}_+$ is called  a
{\it root} if $\alpha \not=0$ and $\geg_\alpha$ is
non-zero. Let $\Delta$ be the set of all roots and
$\Delta_{\pm}=\Delta \cap \pm Q_+$.
For a root $\alpha \in \Delta$ we denote
$\mult_\0o \alpha =\dim \geg_{\alpha,\0o}$,
$\mult_\1o \alpha = -\dim \geg_{\alpha,\1o}$ and
$$
\mult~\alpha=\mult_\0o \alpha +\mult_\1o \alpha=
\dim \geg_{\alpha ,\0o} - \dim \geg_{\alpha,\1o}\ .
\tag{2.5.12}
$$
The set of roots and multiplicities of roots are $W$-invariant. 
We have the formal 
denominator identity for $\geg$ (which is similar to \thetag{2.4.5}): 
$$
\multline
\Phi (z)=\sum_{w \in W}{\epsilon (w)}
\left(
\exp(-2\pi i (w(\rho),z))-
\sum\Sb a \in M \cap \br_{++}\M\endSb 
m(a)\exp(-2\pi i (w(\rho +a), z))\right) =\\ 
=\exp{\left(-2\pi i(\rho,z)\right)}
\prod_{\alpha \in \Delta_+}
{\left( 1-\exp{ \left(-2\pi i (\alpha , z)\right)}\right)^{\mult~\alpha}}. 
\endmultline
\tag{2.5.13}
$$
We remark that $\geg_{m\alpha} =0$ for $m>1$ and  
$\geg_\alpha$ is even of dimensions one if $\alpha$ is a real even root 
(i.e. it belongs to $W(P(\M)_\0o)$). 
Then $\mult~m\alpha =0$ for $m>1$, and 
$\mult~\alpha=1$. If $\alpha$ is a real 
odd root, (i.e. it belongs to $W(P(\M)_\1o)$, the  
space $\geg_{m\alpha}=0$ for $m>2$,  
$\geg_{\alpha}$ is odd of dimension one, and 
$\geg_{2\alpha}$ is even of dimension one. Thus, $\mult~\alpha =-1$,  
$\mult~2\alpha=1$ and $\mult~m\alpha=0$ for $m>2$. 
Suppose that the infinite product \thetag{2.5.13} converges   
at a neighborhood  of the cusp $c$ 
(i.e. if $(\text{Im}~z, \text{Im}~z)\ll 0$). 
It then follows that all zeros of $\Phi$ at a neighborhood of 
the cusp $c$ are quadratic divisors of multiplicity one 
orthogonal to some elements of 
the form $\alpha+mc$, $m \in \bq$, where $\alpha$ is a real root. 

\definition{Definition 2.5.2} A generalized Kac--Moody superalgebra
$\geg^{\prime\prime} ({}_s\Delta)$ defined above 
using the automorphic form $\Phi$ on IV type domain   
with Fourier expansion of Lie type, 
is called an 
{\it automorphic Lorentzian Kac--Moody superalgebra} 
with the symmetrized generalized Cartan matrix 
$B$ (or with a generalized Cartan matrix 
$A$) and with the subset of odd indices $P(\M )_\1o\subset P(\M )$. 
 
The algebra $\geg^{\prime\prime} ({}_s\Delta)$ is also called an   
{\it automorphic correction of a Lorentzian Kac--Moody superalgebra 
$\geg ^{\prime\prime}({}_s\Delta^\re )$, ${}_s\Delta^\re=P(\M )$},  
defined by the symmetrized generalized Cartan matrix 
$B$ (or by a generalized Cartan matrix $A$) and 
with the subset of odd indices $P(\M )_\1o\subset P(\M )$.       
\enddefinition 

To find all possible automorphic generalized 
Lorentzian Kac--Moody superalgebras, one should find all automorphic 
forms $\Phi$ on IV type domains with Fourier expansion of 
Lie type at a cusp. There are two cases: 

\noindent
Case (A): the Weyl group $W$ is non-trivial, equivalently  
$P(\M )\not=\emptyset$. 

\noindent
Case (B): the Weyl group $W$ is trivial, equivalently  
$P(\M )=\emptyset$. 

From Proposition 2.4.3, we get for the case (A): 

\proclaim{Theorem 2.5.3} Assume   
there exists an automorphic form $\Phi$ with 
Fourier expansion of Lie type related with a non-trivial  
reflection group $W$ (or $P(\M)\not=\emptyset$). 
Equivalently,  
suppose that there exists an automorphic correction of the  
Lorentzian Kac--Moody superalgebra 
$\geg^{\prime\prime} (P(\M))$ with $P(\M)\not=\emptyset$.  
Then $P(\M)$ has a lattice Weyl 
vector $\rho \in \br_{++}\M$ and $W$ is elliptic if $(\rho, \rho)<0$, 
and $W$ is parabolic if $(\rho, \rho)=0$. 
\endproclaim

One can consider Theorem 2.5.3 as a generalization to superalgebras 
of some results from \cite{N10} (see also \cite{N9}) 
where it was shown that to have a ``good'' 
Lorentzian Kac--Moody algebra and to have its automorphic correction  
by an automorphic form on IV type domain,   
one is forced to restrict with Weyl groups $W$ of 
elliptic or parabolic type and sets $P(\M)$ with 
a lattice Weyl vector. 
{\it One has no choice} if he wants to restrict with 
holomorphic automorphic forms on IV type domains 
as correcting automorphic forms. 

Type IV domain may have several cusps, and the 
automorphic form $\Phi$ may have different Fourier expansions at  
different cusps.  
Due to R. Borcherds \cite{B5} and 
J. Harvey and G. Moore \cite{HM1}, under some quite general conditions, 
the infinite product \thetag{2.5.13} converges at a 
neighborhood of the cusp and has zeros which are 
quadratic divisors with multiplicity one orthogonal to real roots.  
Globalizing these local data, we suggest a  
definition of ``the most beautiful'' automorphic forms related with 
Kac--Moody algebras: 

\definition{Definition 2.5.4} An automorphic form $\Phi$ on the domain 
$\Omega^+(T)$ where $T$ is a lattice with $2$ negative squares is called a 
{\it Lie reflective automorphic form} if 

\noindent
(i) $\Phi$ is reflective for $T$ and all zeros of 
$\Phi$ have multiplicity one; 

\noindent 
(ii) at any cusp $c$ in $\Omega (T)$, the form 
$\Phi$ has a Lie type Fourier expansion \thetag{2.5.5} with respect to 
a reflection subgroup $W$ of the reflection group of the lattice $T$. 

The lattice $T$  
having a Lie reflective automorphic form is called {\it Lie reflective}.    
\enddefinition 

It is an interesting but very difficult problem to find all Lie reflective  
lattices $T$ and all Lie reflective automorphic forms and the corresponding 
automorphic Kac--Moody superalgebras.  
We believe that because of Theorems 2.2.1, 2.5.3 and Necessary 
Condition 2.2.2 we have the following statement:  

\proclaim{Conjecture 2.5.5} The numbers of Lie reflective 
lattices $T$ with $2$ 
negative squares and their Lie reflective automorphic forms are finite up to 
multiplication of the form of the lattice $T$. 
\endproclaim

This statement is mirror symmetric to finiteness of   
generalized Cartan matrices of elliptic and parabolic type with 
a lattice Weyl vector proved in \cite{N10}.

\smallpagebreak 

We suggest the following physical speculation about  
Definition 2.5.4. The symmetric domain $\Omega^+(T)$ might be related 
with moduli space of some physical theory, a Lie reflective 
automorphic form $\Phi$ might be related 
with some global function (like Lagrangian) which defines this 
physical theory. When we approach to a cusp of $\Omega^+(T)$ 
(that corresponds to a degeneration of the theory), an 
automorphic Lorentzian Kac--Moody algebra defined by Fourier expansion 
of $\Phi$ at the cusp appears. It is related with a quantization  
of the physical theory.

\smallpagebreak   

In Part II, using the methods of \cite{G1}---\cite{G4}, 
we find many Lie reflective automorphic forms and 
Lie reflective lattices which have Lie type Fourier 
expansions with hyperbolic generalized Cartan matrices 
of elliptic and parabolic type classified in Sect. 1. Thus we 
will restrict considering $3$-dimensional case, but all these 
results may be transferred to higher-dimensional case.


\Refs 
\widestnumber\key{vdG2}

\ref
\key Ba 
\by W.L. Baily
\paper Fourier--Jacobi series
\inbook Algebraic groups and discontinuous subgroups.
Proc. Symp. Pure Math. Vol. IX
\eds A. Borel and G.D. Mostow
\publ Amer. Math. Soc.
\publaddr Providence, Rhode Island
\yr 1966
\pages 296--300
\endref


\ref
\key B1 
\by R. Borcherds
\paper Generalized Kac--Moody algebras
\jour J. of Algebra
\vol 115
\yr 1988
\pages 501--512
\endref

\ref
\key B2 
\by R. Borcherds
\paper The monster Lie algebra
\jour Adv. Math.
\vol 83
\yr 1990
\pages 30--47
\endref

\ref
\key B3
\by R. Borcherds
\paper The monstrous moonshine and monstrous Lie superalgebras
\jour Invent. Math.
\vol 109
\yr 1992
\pages 405--444
\endref


\ref
\key B4 
\by R. Borcherds
\paper Sporadic groups and string theory
\inbook Proc. European Congress of Mathem. 1992
\pages 411--421
\endref


\ref
\key B5 
\by R. Borcherds
\paper Automorphic forms on $O_{s+2,2}(\hskip-0.5pt
\br\hskip-0.5pt)$ and
infinite products
\jour Invent. Math. \vol 120
\yr 1995
\pages 161--213
\endref

\ref
\key B6 
\by R. Borcherds
\paper The moduli space of Enriques surfaces and the fake monster Lie
superalgebra
\jour Topology 
\yr 1996 
\vol 35 \issue 3 
\pages 699--710 
\endref

\ref
\key DVV 
\by R. Dijkgraaf, E. Verlinde and H. Verlinde 
\paper Counting dyons in $N=4$ string theory 
\jour hep-th/9607026 
\endref 

\ref
\key EZ 
\by M. Eichler,  D. Zagier
\book The theory of Jacobi forms
\yr 1985
\publ Progress in Math. 55, Birkh\"auser
\endref


\ref
\key FF 
\by A.J. Feingold and I.B. Frenkel
\paper A hyperbolic Kac--Moody algebra and the theory of
Siegel modular forms of genus 2
\jour Math. Ann.
\vol 263
\issue 1
\yr 1983
\pages 87--144
\endref

\ref
\key Fr 
\by E. Freitag
\book Siegelsche Modulfunktionen
\yr 1983
\publ Springer
\endref


\ref
\key vdG1 
\by G. van der Geer
\book Hilbert modular surfaces
\bookinfo Erg. Math. Grenzgeb., 3.Folge, ${\bold 16}$
\yr 1988
\publ Springer Verlag
\endref

\ref
\key vdG2 
\by G. van der Geer
\paper On the geometry of a Siegel modular threefold
\jour Math. Ann.
\vol 260
\yr 1982
\pages 317--350
\endref


\ref\key G1
\by V.A. Gritsenko
\paper Modular forms and moduli spaces of abelian and K3 surfaces
\jour Algebra i Analyz
\vol 6:6
\yr 1994
\pages 65--102
\transl\nofrills  English transl. in
\jour St.Petersburg Math. Jour.
\vol 6:6
\yr 1995
\pages 1179--1208
\endref



\ref\key G2
\by V.A. Gritsenko 
\paper Jacobi functions of n-variables
\jour Zap. Nauk. Sem. LOMI
\vol 168
\yr 1988
\pages 32--45
\lang Russian
\transl\nofrills English transl. in
\jour J\. Soviet Math\.
\vol 53
\yr 1991
\pages 243--252
\endref

\ref
\key G3 
\by V.A. Gritsenko
\paper Arithmetical lifting and its applications
\inbook Number Theory. Proceedings of Paris Seminar  1992--93
\eds S. David
\publ Cambridge Univ. Press
\yr 1995
\pages 103--126
\endref

\ref
\key G4
\by V.A. Gritsenko 
\paper Irrationality of the moduli spaces of polarized abelian
surfaces
\jour The International Mathematics Research Notices
\vol 6
\yr 1994
\pages  235--243,
In  full form  in
``{\it Abelian varieties}'',  Proc. of the  Egloffstein conference
(1993)  de Gruyter, Berlin, 1995, pp. 63--81
\endref


\ref
\key G5
\by V.A. Gritsenko 
\paper Modulformen zur Paramodulgruppe und Modulr\"aume der
Abelschen Variet\"aten
\jour Mathematica Gottingensis Schrift.
des SFB ``Geometrie und Analysis'',
\vol Heft 12
\yr 1995
\pages 1--89
\endref

\ref
\key G6
\by V.A. Gritsenko 
\paper Jacobi functions and Euler products for Hermitian modular
forms
\jour Zap. Nauk. Sem. LOMI
\vol 183 \yr 1990 \pages 77--123
\lang Russian
\transl\nofrills English transl. in
\jour J. Soviet Math.
\yr 1992
\vol 62
\pages 2883--2914
\endref

\ref
\key G7 
\by V.A. Gritsenko 
\book 
Dirichlet series with Euler product in 
the theory of modular forms with respect 
to the orthogonal groups.
\bookinfo Publication of LOMI, E--11--87, 23pp.
\publ ``Nauka"
\publaddr Leningrad
\yr 1987
\endref


\ref\key G8
\by V.A. Gritsenko 
\paper Construction of Hermitian modular forms of genus 2 
from  cusp forms of genus 1
\jour Zap. Nauk. Sem. LOMI
\vol 144
\yr 1985
\pages 51--67
\lang Russian
\transl\nofrills English transl. in
\jour J\. Soviet Math\.
\vol 38
\yr 1987
\pages 2065--2078
\endref

\ref
\key GH 
\by V.A. Gritsenko, K. Hulek
\paper Minimal Siegel modular threefolds
\jour Proceedings of the Cambridge Philosophical Society
\yr 1997 \toappear
\nofrills ; alg-geom/9506017.
\endref

\ref
\key GN1 
\by V.A. Gritsenko  and V.V. Nikulin
\paper Siegel automorphic form correction of some Lorentzi\-an
Kac--Moody Lie algebras
\jour Amer. J. Math.
\yr 1996 \toappear
\nofrills ;  alg-geom/9504006.
\endref

\ref
\key GN2 
\by V.A. Gritsenko  and V.V. Nikulin
\paper Siegel automorphic form correction of a Lorentzian
Kac--Moody algebra
\jour C. R. Acad. Sci. Paris S\'er. A--B
\vol 321
\yr 1995
\pages 1151--1156
\endref

\ref
\key GN3 
\by V.A. Gritsenko  and V.V. Nikulin
\paper K3 surfaces, Lorentzian Kac--Moody algebras and
mirror symmetry
\jour  Math. Res. Lett. \yr 1996 \vol 3 \issue 2 \pages 211--229;
\nofrills  alg-geom/9510008.
 \endref

\ref
\key GN4 
\by V.A. Gritsenko  and V.V. Nikulin
\paper The Igusa modular forms and ``the simplest''
Lorentzian Kac--Moody algebras
\jour Matem. Sbornik 
\yr 1996 \issue 11 \toappear  
\nofrills ; alg-geom/9603010.
\endref

\ref
\key HM1 
\by J. Harvey and G. Moore 
\paper Algebras, BPS-states, and strings 
\jour Nucl. Physics. 
\vol B463 
\yr 1996
\pages 315   
\nofrills hep-th/9510182 
\endref 


\ref
\key HM2  
\by J. Harvey and G. Moore 
\paper On the algebras of BPS-states  
\jour hep-th/9609017   
\yr 1996  
\endref 


\ref
\key Ig1 
\by J. Igusa
\paper On Siegel modular forms of genus two (II)
\jour Amer. J. Math.
\yr 1964
\vol 84
\issue 2
\pages 392--412
\endref

\ref
\key Ig2 
\by J. Igusa
\paper Ring of modular forms of degree two over
$\bz$
\jour Am. J. Math.
\yr 1979
\vol 101
\pages 132--148
\endref

\ref
\key K1 
\by V. Kac
\book Infinite dimensional Lie algebras
\yr 1990
\publ Cambridge Univ. Press
\endref

\ref
\key K2 
\by V. Kac
\paper Lie superalgebras
\jour Adv. Math.
\vol 26
\yr 1977
\pages 8--96
\endref

\ref
\key K3 
\by V. Kac
\paper Infinite-dimensional algebras, Dedekind's $\eta$-function,
classical M\"obius function and the very strange formula
\jour Adv. Math.
\vol 30
\yr 1978
\pages 85--136
\endref

\ref 
\key Ka1 
\by T. Kawai $N=2$ Heterotic string threshold correction, K3 
surfaces and generalized Kac--Moody superalgebra 
\jour Phys. Lett. 
\vol B371
\yr 1996
\page 59   
\nofrills hep-th/9512046  
\endref

\ref 
\key Ka2 
\by T. Kawai 
\paper String duality and modular forms  
\jour Preprint 
\yr 1996 
\nofrills hep-th/9607078    
\endref
 
\ref
\key M1 
\by H. Maass
\paper Die Multiplikatorsysteme zur Siegelschen Modulgruppe
\jour Nachrichten der Akad. der Wissen. G\"ottingen (II.
Math.-Phys.Klasse)
\yr 1964
\vol Nr 11
\pages 125--135
\endref

\ref
\key M2 
\by H. Maass
\paper \"Uber ein Analogon zur Vermutung von Saito-Kurokawa
\jour Invent. math.
\yr 1980
\vol 60
\pages 85--104
\endref

\ref
\key N1 
\by V.V. Nikulin
\paper Finite automorphism groups of K\"ahler K3 surfaces
\jour Trudy Moskov. Mat. Obshch.
\vol 37
\yr 1979 \pages 73--137
\transl\nofrills English transl. in
\jour Trans. Moscow Math. Soc.
\vol 38 \issue 2 \yr 1980
\endref


\ref
\key N2 
\by V.V. Nikulin
\paper Integral symmetric bilinear forms and some of
their geometric applications
\jour Izv. Akad. Nauk SSSR Ser. Mat.
\vol  43
\yr 1979
\pages 111--177
\transl\nofrills English transl. in
\jour Math. USSR Izv.
\vol 14
\yr 1980
\endref

\ref
\key N3 
\by V.V. Nikulin
\paper On the quotient groups of the automorphism groups of
hyperbolic forms by the subgroups generated by 2-reflections,
Algebraic-geometric applications
\jour Current Problems in Math. Vsesoyuz. Inst. Nauchn. i
Tekhn. Informatsii, Moscow
\yr 1981
\pages 3--114
\transl\nofrills English transl. in
\jour J. Soviet Math.
\yr 1983
\vol 22
\pages 1401--1476
\endref

\ref
\key N4 
\by V.V. Nikulin
\paper On arithmetic groups generated by
reflections in Lobachevsky spaces
\jour Izv. Akad. Nauk SSSR Ser. Mat.
\vol  44   \yr 1980 \pages 637--669
\transl\nofrills English transl. in
\jour Math. USSR Izv.
\vol 16 \yr 1981
\endref

\ref
\key N5 
\by V.V. Nikulin
\paper On the classification of arithmetic groups generated by
reflections in Lobachevsky spaces
\jour Izv. Akad. Nauk SSSR Ser. Mat.
\vol  45
\issue 1
\yr 1981
\pages 113--142
\transl\nofrills English transl. in
\jour Math. USSR Izv.
\vol 18
\yr 1982
\endref

\ref
\key N6 
\by V.V. Nikulin
\paper Involutions of integral quadratic forms and their
applications to real algebraic geometry 
\jour Izv. Akad. Nauk SSSR Ser. Mat.
\vol  47 \issue 1  \yr 1983
\transl\nofrills English transl. in
\jour Math. USSR Izv.
\vol 22 \yr 1984 \pages 99--172
\endref

\ref
\key N7 
\by V.V. Nikulin
\paper
Surfaces of type K3 with finite automorphism group and Picard
group of rank three
\jour Proc. Steklov. Math. Inst.
\yr 1984
\vol 165
\pages 113--142
\transl\nofrills English transl. in
\jour Trudy Inst. Steklov
\yr 1985
\vol 3
\endref

\ref
\key N8 
\by V.V. Nikulin
\paper Discrete reflection groups in Lobachevsky spaces and
algebraic surfaces
\inbook Proc. Int. Congr. Math. Berkeley 1986
\vol  1
\pages 654--669
\endref

\ref
\key N9 
\by V.V. Nikulin
\paper A lecture on Kac--Moody Lie algebras of the arithmetic type
\jour Preprint Queen's University, Canada
\vol \#1994-16,
\yr 1994 \nofrills; alg-geom/9412003.
\endref


\ref
\key N10 
\by V.V. Nikulin
\paper Reflection groups in Lobachevsky spaces and
the denominator identity for Lorentzian Kac--Moody algebras
\jour Izv. Akad. Nauk of Russia. Ser. Mat.
\vol  60
\issue 2
\yr 1996
\pages 73--106
\transl\nofrills English transl. in
\jour Russian Acad. Sci. Izv. Math.
\nofrills ; alg-geom/9503003.
\endref

\ref
\key N11 
\by V.V. Nikulin
\paper The remark on discriminants of K3 surfaces moduli as sets
of zeros of automorphic forms \jour Duke e-prints
alg-geom/9512018.
\endref

\ref 
\key N12
\by V.V. Nikulin 
\paper Diagram method for 3-folds and its application to
K\"ahler cone and Picard number of Calabi-Yau 3-folds. I 
\inbook Higher dimensional complex varieties:
Proc. of Intern. Confer. held in Trento, Italy, June 15-24, 1994.
\eds  M. Andreatta, Th. Peternell 
\publ de Gruyter 
\yr 1996 
\pages 261--328
\endref 


\ref
\key R
\by U. Ray
\paper A character formula for generalized Kac--Moody superalgebras
\jour J. of Algebra
\vol 177
\yr 1995
\pages 154--163
\endref

\ref
\key V1 
\by \'E.B. Vinberg 
\paper On groups of unit elements of certain quadratic forms
\jour Mat. Sbornik
\yr 1972
\vol 87
\pages 18--36
\transl\nofrills English transl. in
\jour Math USSR Sbornik
\vol 16
\yr 1972
\pages 17--35
\endref

\ref
\key V2 
\by \'E.B. Vinberg 
\paper The absence of crystallographic reflection groups in Lobachevsky
spaces of large dimension
\jour Trudy Moscow. Mat. Obshch.
\vol  47 \yr 1984  \pages 68 -- 102
\transl\nofrills English transl. in
\jour Trans. Moscow Math. Soc.
\vol 47 \yr 1985
\endref


\ref
\key V3 
\by \'E.B. Vinberg 
\paper Hyperbolic reflection groups
\jour Uspekhi Mat. Nauk
\vol 40
\yr 1985
\pages 29--66
\transl\nofrills English transl. in
\jour Russian Math. Surveys
\vol 40
\yr 1985
\endref



\endRefs
\enddocument

\end